\documentclass{pspum-l}
\usepackage{amscd,amssymb,graphics,color,hyperref}

\usepackage{mathrsfs}

\usepackage[all]{xy}

\newtheorem{theorem}{Theorem}[section]
\newtheorem{lemma}[theorem]{Lemma}
\newtheorem{proposition}[theorem]{Proposition}

\newtheorem{conjecture}[theorem]{Conjecture}
\newtheorem{xca}[theorem]{Exercise}
\newtheorem{question}[theorem]{Question}

\theoremstyle{definition}
\newtheorem{definition}[theorem]{Definition}

\newtheorem{example}[theorem]{Example}
\newtheorem{examples}[theorem]{Examples}

\theoremstyle{remark}
\newtheorem{remark}[theorem]{Remark}

\newcommand{\NN} {\mathbb{N}}
\newcommand{\ZZ} {\mathbb{Z}}
\newcommand{\QQ} {\mathbb{Q}}
\newcommand{\RR} {\mathbb{R}}
\newcommand{\CC} {\mathbb{C}}

\newcommand{\PP} {\mathbb{P}}
\renewcommand{\AA} {\mathbb{A}}

\newcommand {\shC} {\mathcal{C}}

\newcommand {\shExt} {\mathcal{E} \!\text{\textit{xt}}}
\newcommand {\shE} {\mathcal{E}}

\newcommand {\shL} {\mathcal{L}}
\newcommand {\shM} {\mathcal{M}}

\newcommand {\shN} {\mathcal{N}}

\newcommand {\shT} {\mathcal{T}}

\newcommand {\shP} {\mathcal{P}}

\newcommand {\shX} {\mathcal{X}}

\renewcommand {\Bar} {\operatorname{Bar}}

\newcommand {\Di} {{\rm D}}

\newcommand {\Dlog} {\operatorname{Dlog}}

\newcommand {\dual} {\vee}

\newcommand {\gp} {{\operatorname{gp}}}

\newcommand {\Hom} {\operatorname{Hom}}

\renewcommand {\Im} {\operatorname{Im}}

\newcommand {\Int} {\operatorname{Int}}

\newcommand {\M} {\mathcal{M}}
\newcommand {\Mbar} {\overline{\M}}

\renewcommand{\O} {\mathcal{O}}

\renewcommand{\P} {\mathscr{P}}

\newcommand {\Pic} {\operatorname{Pic}}

\newcommand {\pre} {\mathrm{pre}}

\newcommand {\Proj} {\operatorname{Proj}}

\newcommand {\rank} {\operatorname{rank}}

\newcommand {\Spec} {\operatorname{Spec}}

\newcommand {\T} {\shT}
\newcommand {\X} {\shX}

\newcommand {\U} {\mathscr{U}}

\newcommand {\shLS} {\mathcal{LS}}

\def\mydate{\ifcase\month \or January\or February\or March\or
April\or May\or June\or July\or August\or September\or October\or 
November\or December\fi \space\number\day,\space\number\year}

\begin{document}
\def\mapright#1{\smash{
 \mathop{\longrightarrow}\limits^{#1}}}
\def\mapleft#1{\smash{
 \mathop{\longleftarrow}\limits^{#1}}}
\def\exact#1#2#3{0\to#1\to#2\to#3\to0}
\def\mapup#1{\Big\uparrow
  \rlap{$\vcenter{\hbox{$\scriptstyle#1$}}$}}
\def\mapdown#1{\Big\downarrow
  \rlap{$\vcenter{\hbox{$\scriptstyle#1$}}$}}
\def\dual#1{{#1}^{\scriptscriptstyle \vee}}
\def\invlim{\mathop{\rm lim}\limits_{\longleftarrow}}
\def\cy{\check y}

\title[The Strominger-Yau-Zaslow conjecture.]{
The Strominger-Yau-Zaslow conjecture: From torus fibrations to
degenerations.}

\author{Mark Gross} 
\address{UCSD Mathematics, 9500 Gilman Drive, La Jolla, CA 92093-0112, USA}
\email{mgross@math.ucsd.edu}
\thanks{This work was partially supported by NSF grant
0505325}

\date{}
\begin{abstract}
We trace progress and thinking about the Strominger-Yau-Zaslow conjecture
since its introduction in 1996. In particular, we aim to explain how
the conjecture led to the algebro-geometric program developed by myself
and Siebert, whose objective is
to explain mirror symmetry by studying degenerations
of Calabi-Yau manifolds. We end by outlining how tropical curves arise
in the mirror symmetry story.
\end{abstract}
\maketitle
\bigskip

\section*{Introduction.}
Up to the summer of 1996, there had been a number of spectacular
successes in mirror symmetry. After the pioneering
initial work of Candelas,
de la Ossa, Greene and Parkes \cite{COGP} which calculated the instanton
predictions for the quintic three-fold, Batyrev \cite{Bat} gave a powerful
mirror symmetry construction for Calabi-Yau hypersurfaces in toric
varieties, later generalized by Batyrev and Borisov \cite{BB} to complete
intersections in toric varieties. Kontsevich \cite{Kstable} introduced his
notion of stable maps of curves
and studied their moduli, setting the stage
for a new flowering of enumerative geometry. Eventually, this led
to the mathematical calculation of Gromov-Witten
invariants for the quintic in varying forms, 
\cite{Givental,LLY,Bert,Gath}. In between, a great
deal of the structure of mirror symmetry was elucidated by many
researchers in both string theory and algebraic geometry.

On the other hand, at that time
I had been primarily interested in the geometry of Calabi-Yau manifolds.
Many of the results about mirror symmetry seemed to rely primarily
on information about the ambient toric varieties; in particular,
Givental's work \cite{Givental}
and succeeding work by Lian, Liu and Yau \cite{LLY},
Bertram \cite{Bert} and Gathmann \cite{Gath}
always performed calculations on the moduli space
of stable maps into $\PP^4$ rather than the quintic, and in the 
Batyrev-Borisov type constructions, while there was a clear
combinatorial relationship
between the ambient toric varieties, there was no apparent 
geometric relationship
between the Calabi-Yaus themselves. As a result, I tended to avoid
thinking about mirror symmetry precisely because of this lack of
geometric understanding.

This changed dramatically in June of 1996, when Strominger, Yau
and Zaslow \cite{SYZ} released their paper ``Mirror Symmetry is $T$-duality.''
They made a remarkable proposal, based on recent ideas in string theory,
that for the first time gave a geometric interpretation for mirror
symmetry. 

Let me summarize, very roughly, the physical argument here. Developments
in string theory in the mid-1990s had introduced the notion of \emph{Dirichlet
branes}, or $D$-branes. These are submanifolds of space-time, with
some additional data, which
should serve as a boundary condition for open strings, i.e.\ we allow
open strings to propagate with their endpoints constrained to lie on
a $D$-brane. Remembering that space-time, according to string theory, 
looks like
$\RR^{1,3}\times X$, where $\RR^{1,3}$ is ordinary space-time and
$X$ is a Calabi-Yau three-fold, we can split a $D$-brane into a product
of a submanifold of $\RR^{1,3}$ and one on $X$. It turned out, simplifying
a great deal, that there were 
two particular types of submanifolds on $X$ of interest: 
\emph{holomorphic} $D$-branes,
i.e.\ holomorphic submanifolds with a holomorphic line bundle, and 
\emph{special Lagrangian} $D$-branes, which are 
\emph{special Lagrangian submanifolds} with 
flat $U(1)$ bundle: 

\begin{definition} Let $X$ be an $n$-dimensional
Calabi-Yau manifold with $\omega$
the K\"ahler form of a Ricci-flat metric on $X$ and $\Omega$
a nowhere vanishing holomorphic $n$-form. Then a submanifold $M\subseteq
X$ is \emph{special Lagrangian} if it is Lagrangian, i.e.\ 
$\dim_{\RR} M=\dim_{\CC} X$
and $\omega|_M=0$, and in addition $\Im\Omega|_M=0$.
\end{definition}

The origins of mirror symmetry in physics suggest that if $X$ and $\check X$
are a mirror pair of Calabi-Yau manifolds, then string theory on a 
compactification of space-time using $X$ should be the same as that using
$\check X$, but with certain data interchanged. In the case of $D$-branes,
the suggestion is that the moduli space of holomorphic $D$-branes on $X$
should be isomorphic to the moduli space of special Lagrangian $D$-branes 
on $\check X$.
Now $X$ itself is the moduli space of points on $X$. So each
point on $X$ should correspond to a pair $(M,\nabla)$, where $M\subseteq\check
X$ is a special Lagrangian submanifold and $\nabla$ is a flat connection on
$M$. 

A theorem of McLean \cite{McLean}
tells us that the tangent space to the moduli space of
special Lagrangian deformations of a special Lagrangian submanifold $M\subseteq
\check X$ is $H^1(M,\RR)$. Of course, 
the moduli space of flat $U(1)$-connections
modulo gauge equivalence on $M$ is the torus $H^1(M,\RR)/H^1(M,\ZZ)$.
In order for this moduli space to be of the correct dimension, we need
$\dim H^1(M,\RR)=n$, the complex dimension of $X$. This suggests that 
$X$ consists of a family of tori which are dual to a family of special 
Lagrangian tori on $\check X$. An elaboration of this argument yields the
following conjecture:

\begin{conjecture}
\emph{The Strominger-Yau-Zaslow conjecture}. If $X$ and $\check X$ are
a mirror pair of Calabi-Yau $n$-folds, then there exists fibrations
$f:X\rightarrow B$ and $\check f:\check X\rightarrow B$ whose fibres
are special Lagrangian, with general fibre an $n$-torus. Furthermore,
these fibrations are dual, in the sense that canonically
$X_b=H^1(\check X_b,\RR/\ZZ)$
and $\check X_b=H^1(X_b,\RR/\ZZ)$ whenever $X_b$ and $\check X_b$ 
are non-singular tori.
\end{conjecture}

I will clarify this statement as we review the work of the past
ten years; however, as I have stated this conjecture,
it is likely to be false.  On the other hand, there are weaker
versions of the conjecture which probably are true. Even better,
these weaker statements are
probably within reach of modern-day technology (with a lot
of hard work). Nevertheless, there has been a lot of good progress
on precise versions of the above conjecture at the topological and
symplectic level. In addition, the conjecture has been successful
at explaining many features of mirror symmetry, some of which are
still heuristic and some of which are rigorous. It is my belief that
a final satisfactory understanding of mirror symmetry will flow from
the SYZ conjecture, even if results do not take the form initially suggested
by it.

My main goal here is to explain the journey taken over the last ten years.
I want to focus on explaining the evolution and development of the
ideas, rather than focus on precise statements. Except in the first
few sections, I will give few precise statements.

In those first sections, I will clarify
the above statement of the conjecture, and show how
it gives a satisfactory explanation of mirror symmetry in the so-called
semi-flat case, i.e.\ the case when the metric along the special Lagrangian
fibres is flat. This leads naturally to a discussion of affine manifolds,
metrics on them, and the Legendre transform. These now appear to be
the key structures underlying mirror symmetry. 

We next take a look at the case when singular fibres appear. In this case
we need to abandon the precise form of duality we developed in the semi-flat
case and restrict our attention to topological duality. In the realm
of purely topological duality, the SYZ conjecture has been entirely successful
at explaining topological features of mirror symmetry for a large range
of Calabi-Yau manifolds, including those produced by the Batyrev-Borisov
construction for complete intersections in toric varieties.

Moving on, we take a look at Dominic Joyce's arguments demonstrating the
problems with the strong form of the SYZ conjecture stated above. This forces
us to recast the SYZ conjecture as a limiting statement. Mirror symmetry
is always about the behaviour of Calabi-Yau manifolds near maximally unipotent
degenerations. A limiting form of the SYZ conjecture
suggests that one can find special Lagrangian tori on Calabi-Yau
manifolds near a maximally unipotent degeneration, and as we approach
the limit point in complex moduli space, we expect to see a larger portion
of the Calabi-Yau manifold filled out by special Lagrangian tori. Unlike
the original SYZ conjecture, though still difficult,
this one looks likely to be accessible by current techniques.

This form of the conjecture then motivates a new round of questions.
In this limiting picture, we expect the base $B$ of the hypothetical
special Lagrangian fibration to be the so-called Gromov-Hausdorff limit
of a sequence of Calabi-Yau manifolds approaching the maximally unipotent
degeneration. Gromov-Hausdorff convergence is a metric space concept,
while maximally unipotent degeneration is an algebro-geometric, Hodge-theoretic
concept. How do these two concepts relate? In the summer of 2000,
Kontsevich suggested that the Gromov-Hausdorff limit will be, roughly,
the dual intersection complex of the algebro-geometric degeneration, at
least on a topological level. We explore this idea in \S 6.

On the other hand, how does this help us with mirror symmetry? Parallel
to these developments on limiting forms of the SYZ conjecture, my coauthor
Bernd Siebert had been studying degenerations of Calabi-Yau manifolds
using logarithmic geometry with Stefan Schr\"oer. Siebert noticed that
mirror symmetry seemed to coincide with a combinatorial exchange of logarithmic
data on the one side and polarizations on the other. Together, we realised
that this approach to mirror symmetry meshed well with the limiting picture
predicted by SYZ. Synthesizing these two approaches, we discovered
an algebro-geometric version of the SYZ approach, which I will describe
here. The basic idea is to forget about special Lagrangian fibrations,
and only keep track of the base of the fibration, which is an affine
manifold. We show how polyhedral decompositions of affine manifolds give
rise to degenerate Calabi-Yau varieties, and conversely how certain
sorts of degenerations of Calabi-Yau varieties, which we call \emph{toric
degenerations}, give rise to affine manifolds as their dual intersection
complex. Mirror symmetry is
again explained by a discrete version of the  Legendre transform much 
as in the semi-flat case.

We end with a discussion of the connection of tropical curves with
this approach. The use of tropical curves in curve counting in two-dimensional
toric varieties has been pioneered
in work of Mihkalkin \cite{Mik}; Nishinou and Siebert \cite{NS}
generalized
this work to higher dimensions using an approach directly inspired by
the approach I discuss here. On the other hand, tropical curves have not
yet been used for counting curves in Calabi-Yau manifolds, so I will end
the paper by discussing how tropical curves arise naturally in our picture.
\medskip

I would like to thank the organizers of the Seattle conference for
running an excellent conference, and my coauthors Pelham Wilson and
Bernd Siebert on SYZ related results; much of the work mentioned here 
came out of work with them.
\bigskip

\section{First a topological observation.}

Before doing anything else, let's ask a very basic question: why should
dualizing torus fibrations interchange Hodge numbers of Calabi-Yau threefolds?
If you haven't seen this, it's the first thing one should look at as it
is particularly easy to see, if we make a few assumptions. Suppose we
are given a pair of Calabi-Yau threefolds $X$ and $\check X$ with
fibrations $f:X\rightarrow B$, $\check f:\check X\rightarrow B$ with
the property that there is a dense open set $B_0\subseteq B$ such that
$f_0:f^{-1}(B_0)\rightarrow B_0$ and $\check f_0:\check f^{-1}(B_0)\rightarrow
B_0$ are torus fibre bundles. So all the singular fibres of $f$ and $\check f$
lie over $\Gamma:=B\setminus B_0$. (Note that unless $\chi(X)=0$, there
must be some singular fibres.) Finally, assume $f_0$ and $\check f_0$
are dual torus fibrations, i.e.\ $f_0$ can be identified with the torus
fibration $R^1\check f_{0*}(\RR/\ZZ)\rightarrow B_0$ and $\check f_0$
can be identified with the torus fibration $R^1 f_{0*}(\RR/\ZZ)\rightarrow B_0$.
(This is a slight abuse of notation: by $R^1\check f_{0*}(\RR/\ZZ)$
we really mean the torus bundle
obtained by taking the vector bundle associated to the local system
$R^1\check f_{0*}\RR$ and dividing out by the family of lattices 
$R^1\check f_{0*}\ZZ$.)

If $V/\Lambda$ is a single torus with $V$ an $n$-dimensional
vector space and $\Lambda$
a lattice in $V$, then $H^p(V/\Lambda,\RR)\cong \bigwedge^p \dual{V}$, while the
dual torus, $\dual{V}/\dual{\Lambda}$, has $H^p(\dual{V}/\dual{\Lambda},
\RR)\cong \bigwedge^p V$. If we choose an isomorphism $\bigwedge^n V\cong
\RR$, then we get an isomorphism $H^p(V/\Lambda,\RR)\cong H^{n-p}(\dual{V}/
\dual{\Lambda},\RR)$. Similarly, in the relative setting for $f_0$ and
$\check f_0$, if we have an isomorphism $R^3f_{0*}\RR\cong \RR$, we obtain
isomorphisms
\[
R^pf_{0*}\RR\cong R^{3-p}\check f_{0*}\RR.
\]
We now make a simplifying assumption. Let $i:B_0\hookrightarrow B$
be the inclusion. We will say $f$ is \emph{$\RR$-simple} if
\[
i_*R^pf_{0*}\RR\cong R^pf_*\RR
\]
for all $p$. (We can in general replace $\RR$ by any abelian group
$G$, and then we say $f$ is $G$-simple.) Of course, not all torus
fibrations are $\RR$-simple, but it turns out that the most interesting
ones which occur in the topological form of SYZ are. So let's 
assume $f$ and $\check f$ are $\RR$-simple. 
With this assumption, we obtain isomorphisms
\begin{equation}
\label{dualityiso}
R^pf_*\RR\cong R^{3-p}\check f_*\RR.
\end{equation}
We can now use this to study the Leray spectral sequence for $f$
and $\check f$. 

Let's make an additional assumption that $X$ and $\check X$ are
simply connected. So in particular $B$ is simply connected. Let's
assume $B$ is a three-manifold. So we have the $E_2$
terms in the Leray spectral sequence for $f$:
\[
\begin{matrix}
\RR&0&0&\RR\\
H^0(B,R^2f_*\RR)&H^1(B,R^2f_*\RR)&H^2(B,R^2f_*\RR)&H^3(B,R^2f_*\RR)\\
H^0(B,R^1f_*\RR)&H^1(B,R^1f_*\RR)&H^2(B,R^1f_*\RR)&H^3(B,R^1f_*\RR)\\
\RR&0&0&\RR
\end{matrix}
\]
Since $X$ is simply connected, $H^1(X,\RR)=H^5(X,\RR)=0$, from which
we conclude that $H^0(B,R^1f_*\RR)=H^3(B,R^2f_*\RR)=0$. The same
argument works for $\check f$, and then (\ref{dualityiso}) gives
$H^0(B,R^2f_*\RR)\cong H^0(B,R^1\check f_*\RR)=0$ and similarly
$H^3(B,R^1f_*\RR)=0$. Finally, consider the possible non-zero maps
for the spectral sequence:
\[
\xymatrix@C=30pt
{\RR\ar[rrd]^{d_1}&0&0&\RR\\
0&H^1(B,R^2f_*\RR)&H^2(B,R^2f_*\RR)&0\\
0&H^1(B,R^1f_*\RR)\ar[rrd]^{d_2}&H^2(B,R^1f_*\RR)&0\\
\RR&0&0&\RR}
\]
We need one further assumption, which is again natural given the duality
relationship of $f_0$ and $\check f_0$: both $f$ and $\check f$ possess
sections. (Actually, working over $\RR$, we just need the existence
of cohomology classes on $X$ and $\check X$
which evaluate to something non-zero on a fibre of $f$ or $\check f$.)
A section intersects each fibre non-trivially, and hence gives 
a section of $R^3f_*\RR$. Since such a section also represents a cohomology
class on $X$, $d_1$ must be the zero map. Similarly, a fibre of $f$
cannot be homologically trivial because it intersects the section non-trivially,
and thus the map $d_2$ must be zero. So the spectral sequence degenerates
at $E_2$. In particular, we get $\RR^{h^{1,1}}\cong
H^2(X,\RR)\cong H^1(B,R^1f_*\RR)\cong H^1(B,R^2\check f_*\RR)$ and 
$\RR^{h^{2,2}}\cong H^4(X,\RR)\cong H^2(B,R^2f_*\RR)\cong H^2(B,R^1\check
f_*\RR)$, where $h^{p,q}$ are the Hodge numbers of $X$. Thus the
third Betti number of $\check X$ is $2+h^{1,1}+h^{2,2}=2(1+h^{1,1})$,
so we see $h^{1,1}(X)=h^{1,2}(\check X)$ and $h^{1,2}(X)
=h^{1,1}(\check X)$. 

So modulo some assumptions which would of course eventually have to be
justified, it is clear, at least in the three-dimensional case, why the
Hodge numbers are interchanged by duality. 

More generally, in any dimension, one might hope that $\RR$-simplicity
implies $\dim_{\RR} H^p(B,R^qf_*\RR)=h^{p,q}$, and then a more general
exchange of Hodge numbers becomes clear. See Theorem \ref{hodgedecomp}
for a related result.

\begin{remark}
This argument can be refined over $\ZZ$ to make new predictions about
the behaviour of \emph{integral} cohomology under mirror symmetry.
In \cite{SlagII}, Theorem 3.10, it was shown, again in the
three-dimensional case, that if $f$ and $\check f$ are
$\ZZ$-simple and $\QQ/\ZZ$-simple, $f$ and $\check f$ have sections,
and $H^1(X,\ZZ)=0$, then
\begin{eqnarray*}
H^{even}(X,\ZZ[1/2])&\cong&H^{odd}(\check X,\ZZ[1/2])\\
H^{odd}(X,\ZZ[1/2])&\cong&H^{even}(\check X,\ZZ[1/2]).
\end{eqnarray*}
There are problems in the argument with two-torsion, but it is likely
the above isomorphisms hold over $\ZZ$. See \cite{BS}
for evidence for this latter conjecture.
\end{remark}

Enough speculation.
Now let's get serious about the structure of special Lagrangian fibrations.

\section{Moduli of special Lagrangian submanifolds}

The first step in really understanding the SYZ conjecture is to 
examine the structures which arise on the base of a special Lagrangian
fibration. These structures arise from McLean's theorem on the moduli
space of special Lagrangian submanifolds \cite{McLean}, and these
structures and their relationships were explained by Hitchin in \cite{Hit}. 
We outline some of these ideas here.
McLean's theorem says that the moduli space of deformations of a compact
special Lagrangian submanifold
of a compact 
Calabi-Yau manifold $X$ is unobstructed, with tangent space at $M\subseteq
X$ special Lagrangian canonically isomorphic to
the space of harmonic $1$-forms on $M$. This isomorphism is seen
explicitly as follows. Let $\nu\in\Gamma(M,N_{M/X})$ be a normal
vector field to $M$ in $X$. Then $(\iota(\nu)\omega)|_M$ and
$(\iota(\nu)\Im\Omega)|_M$ are both seen to be well-defined forms
on $M$: one needs to lift $\nu$ to a vector field but the choice is
irrelevant because $\omega$ and $\Im\Omega$ restrict to zero on $M$.
McLean shows that if $M$ is special Lagrangian then
\[
\iota(\nu)\Im\Omega=-*\iota(\nu)\omega,
\]
where $*$ denotes the Hodge star operator on $M$, and furthermore,
$\nu$ corresponds to an infinitesimal deformation preserving the
special Lagrangian condition if and only if $d(\iota(\nu)\omega)
=d(\iota(\nu)\Im\Omega)=0$. This gives the correspondence between
harmonic $1$-forms and infinitesimal special Lagrangian deformations.

Let $f:X\rightarrow B$ be a special Lagrangian fibration
with torus fibres, and assume for now that all fibres of $f$ are non-singular.
Then we obtain three structures on $B$: two affine structures and a metric,
as we shall now see.

\begin{definition} 
\label{affine}
Let $B$ be an $n$-dimensional manifold.
An {\it affine structure} on $B$ is given by an atlas $\{(U_i,\psi_i)\}$
of coordinate charts $\psi_i:U_i\rightarrow \RR^n$,
whose transition functions $\psi_i\circ\psi_j^{-1}$ lie in ${\rm Aff}(\RR^n)$.
We say the affine structure is \emph{tropical} if the transition functions
lie in $\RR^n\rtimes GL(\ZZ^n)$, i.e.\ have integral linear part. We say
the affine structure is {\it integral} if the transition functions
lie in ${\rm Aff}(\ZZ^n)$. 

If an affine manifold $B$ carries a Riemannian metric $g$, then we say
the metric is \emph{affine K\"ahler} or \emph{Hessian} if $g$ is locally
given by $g_{ij}=\partial^2K/\partial y_i\partial y_j$ for some convex
function $K$ and $y_1,\ldots,y_n$ affine coordinates.
\end{definition}

Then we obtain the three structures as follows:

\emph{Affine structure 1.} For a normal vector field $\nu$ to a fibre $X_b$
of $f$, $(\iota(\nu)\omega)|_{X_b}$ is a well-defined $1$-form on $X_b$,
and we can compute its periods as follows.
Let $U\subseteq B$ be a small open set,
and suppose we have submanifolds $\gamma_1,\ldots,\gamma_n\subseteq
f^{-1}(U)$ which are families of 1-cycles over $U$ and such that
$\gamma_1\cap X_b,\ldots,\gamma_n\cap X_b$ 
form a basis for $H_1(X_b,\ZZ)$ for each
$b\in U$. Consider the $1$-forms $\omega_1,\ldots,\omega_n$ on $U$
defined by fibrewise integration:
\[
\omega_i(\nu)=\int_{X_b\cap\gamma_i} \iota(\nu)\omega,
\]
for $\nu$ a tangent vector on $B$ at $b$, which we can lift
to a normal vector field of $X_b$. We have $\omega_i=f_*(\omega|_{\gamma_i})$,
and since $\omega$ is closed, so is $\omega_i$. Thus there are locally
defined functions $y_1,\ldots,y_n$ on $U$ with $dy_i=\omega_i$. 
Furthermore, these functions are well-defined up to the choice of basis
of $H_1(X_b,\ZZ)$ and constants. Finally, they give well-defined coordinates,
as follows from the fact that 
$\nu\mapsto \iota(\nu)\omega$ yields an isomorphism of $\T_{B,b}$ with
$H^1(X_b,\RR)$ by McLean's theorem. Thus $y_1,\ldots,y_n$ define local
coordinates of a tropical affine structure on $B$.

\emph{Affine structure 2.} We can play the same trick with $\Im\Omega$: 
choose submanifolds $\Gamma_1,\ldots,\Gamma_n\subseteq f^{-1}(U)$
which are families of $n-1$-cycles over $U$ and such that
$\Gamma_1\cap X_b,\ldots,\Gamma_n\cap X_b$ form a basis for $H^{n-1}(X_b,
\ZZ)$. We define $\lambda_i$ by $\lambda_i=-f_*(\Im\Omega|_{\Gamma_i})$,
or equivalently,
\[
\lambda_i(\nu)=-\int_{X_b\cap\Gamma_i} \iota(\nu)\Im\Omega.
\]
Again $\lambda_1,\ldots,\lambda_n$ are closed $1$-forms, with
$\lambda_i=d\check y_i$ locally, and again $\check y_1,\ldots,\check
y_n$ are affine coordinates for a tropical affine structure on $B$.

\emph{The McLean metric.} The Hodge metric on $H^1(X_b,\RR)$ is given
by 
\[
g(\alpha,\beta)=\int_{X_b} \alpha\wedge *\beta
\]
for $\alpha$, $\beta$ harmonic $1$-forms, and hence induces a
metric on $B$, which can be written as 
\[
g(\nu_1,\nu_2)=-\int_{X_b}\iota(\nu_1)\omega\wedge \iota(\nu_2)\Im\Omega.
\]

\medskip

A crucial observation of Hitchin \cite{Hit} is that these structures are
related by the Legendre transform:

\begin{proposition}
\label{hessianmetric}
Let $y_1,\ldots,y_n$ be local affine coordinates
on $B$ with respect to the affine structure induced by $\omega$. 
Then locally there is a function $K$ on $B$
such that
\[
g(\partial/\partial y_i,\partial/\partial y_j)=\partial^2 K/\partial y_i
\partial y_j.
\]
Furthermore, $\cy_i=\partial K/\partial y_i$
form a system of affine coordinates with respect to the affine
structure induced by $\Im\Omega$, and if
\[
\check K(\cy_1,\ldots,\cy_n)=\sum \cy_i y_i-K(y_1,\ldots,y_n)
\]
is the Legendre transform of $K$, then
\[
y_i=\partial \check K/\partial\cy_i
\]
 and
\[
\partial^2\check K/\partial y_i\partial y_j=g(\partial/\partial\cy_i,
\partial/\partial\cy_j).
\]
\end{proposition}

\begin{proof}
Take families $\gamma_1,\ldots,\gamma_n,\Gamma_1,\ldots,\Gamma_n$
as above
over an open neighbourhood $U$ with the two bases being
Poincar\'e dual, i.e.\ $(\gamma_i\cap X_b)\cdot(\Gamma_j\cap X_b)=
\delta_{ij}$ for $b\in U$.
Let $\gamma_1^*,\ldots,\gamma_n^*$ and $\Gamma_1^*,\ldots,\Gamma_n^*$
be the dual bases for $\Gamma(U,R^1f_*\ZZ)$ and $\Gamma(U,R^{n-1}f_*\ZZ)$
respectively. From the choice of $\gamma_i$'s, we get local coordinates
$y_1,\ldots,y_n$ with $dy_i=\omega_i$, so in particular
\[
\delta_{ij}=\omega_i(\partial/\partial y_j)=\int_{\gamma_i\cap X_b}
\iota(\partial/\partial y_j)\omega,
\]
so $\iota(\partial/\partial y_j)\omega$ defines the cohomology
class $\gamma_j^*$ in $H^1(X_b,\RR)$. Similarly,
let
\[
g_{ij}=-\int_{\Gamma_i\cap X_b}\iota(\partial/\partial y_j)\Im
\Omega;
\]
then $-\iota(\partial/\partial y_j)\Im\Omega$ defines the cohomology
class $\sum_i g_{ij}\Gamma_i^*$ in $H^{n-1}(X_b,\RR)$, and
$\lambda_i=\sum_j g_{ij}dy_j$.
Thus
\begin{eqnarray*}
g(\partial/\partial y_j,\partial/\partial y_k)&=&
-\int_{X_b}\iota(\partial/\partial y_j)\omega
\wedge \iota(\partial/\partial y_k)\Im\Omega\\
&=&g_{jk}.
\end{eqnarray*}
On the other hand, let $\cy_1,\ldots,\cy_n$ be coordinates with
$d\cy_i=\lambda_i$. Then
\[
{\partial\cy_i/\partial y_j}=g_{ij}=g_{ji}={\partial\cy_j/
\partial y_i},
\]
so $\sum\cy_i dy_i$ is a closed 1-form. Thus there exists
locally a function $K$ such that $\partial K/\partial y_i=\cy_i$ and
$\partial^2 K/\partial y_i\partial y_j=g(\partial/\partial y_i,
\partial/\partial y_j)$.
A simple calculation then confirms that $\partial\check K/\partial \cy_i
=y_i$. On the other hand,
\begin{eqnarray*}
g(\partial/\partial\cy_i,\partial/\partial\cy_j)&=&
g\left(\sum_k {\partial y_k\over\partial\cy_i}{\partial\over
\partial y_k},\sum_l {\partial y_l\over\partial\cy_j}
{\partial\over\partial y_l}\right)\\
&=&\sum_{k,l}{\partial y_k\over\partial\cy_i}{\partial y_l\over\partial
\cy_j} g(\partial/\partial y_k,\partial/\partial y_l)\\
&=&\sum_{k,l} {\partial y_k\over\partial\cy_i}{\partial y_l\over
\partial\cy_j}{\partial\cy_k\over\partial y_l}\\
&=&{\partial y_j\over\partial\cy_i}={\partial^2\check K\over
\partial\cy_i\partial\cy_j}.
\end{eqnarray*}
\end{proof}

Thus we introduce the notion of \emph{Legendre transform} of an
affine manifold with a multi-valued convex function.

\begin{definition}
\label{multivaluedconvex}
Let $B$ be an affine manifold. A \emph{multi-valued} function $K$
on $B$ is a collection of functions on an open cover $\{(U_i,K_i)\}$
such that on $U_i\cap U_j$, $K_i-K_j$ is affine linear. We say $K$
is \emph{convex} if the Hessian
$(\partial^2 K_i/\partial y_j\partial y_k)$ is positive definite for
all $i$, in any, or equivalently all, affine coordinate systems $y_1,
\ldots,y_n$.

Given a pair $(B,K)$ of affine manifold and convex multi-valued
function, the \emph{Legendre transform} of $(B,K)$ is a pair $(\check B,
\check K)$ where $\check B$ is an affine structure on the underlying
manifold of $B$ with coordinates given locally
by $\check y_i=\partial K/\partial y_i$, and $\check K$ is defined
by
\[
\check K_i(\check y_1,\ldots,\check y_n)=\sum \check y_j y_j
-K_i(y_1,\ldots,y_n).
\]
\end{definition}

\begin{xca}
Check that $\check K$ is also convex, and that the Legendre transform
of $(\check B,\check K)$ is $(B,K)$.
\end{xca}

\section{Semi-flat mirror symmetry}

Now let's forget about special Lagrangian fibrations for the
moment. Instead, we see how the structures found on $B$ give a toy
version of mirror symmetry.

\begin{definition} Let $B$ be a tropical affine manifold.
\begin{enumerate}
\item
Define $\Lambda\subseteq\T_B$ to be the local system of lattices
generated locally by $\partial/\partial y_1,\ldots,\partial/\partial y_n$, where
$y_1,\ldots,y_n$ are local affine coordinates. This is well-defined because
transition maps are in $\RR^n\rtimes GL_n(\ZZ)$. Set 
\[
X(B):=\T_B/\Lambda;
\]
this is a torus bundle over $B$. In addition, $X(B)$ carries a complex
structure defined locally as follows. Let $U\subseteq B$ be
an open set with affine coordinates $y_1,\ldots,y_n$, so $\T_U$
has coordinate functions $y_1,\ldots,y_n$, $x_1=dy_1,\ldots,x_n=dy_n$.
Then 
\[
q_j=e^{2\pi i(x_j+iy_j)}
\]
gives a system of holomorphic coordinates on $T_U/\Lambda|_U$, and
the induced complex structure is 
independent of the choice of affine coordinates.

Later we will need a variant of this: for $\epsilon>0$, set
\[
X_{\epsilon}(B):=\T_B/\epsilon\Lambda;
\]
this has a complex structure with coordinates given by
\[
q_j=e^{2\pi i(x_j+iy_j)/\epsilon}.
\]
(As we shall see later, the limit $\epsilon\rightarrow 0$ 
corresponds to a large complex structure limit.)
\item
Define $\check\Lambda\subseteq\T^*_B$ to be the local system of
lattices generated locally by $dy_1,\ldots,dy_n$, with $y_1,\ldots,y_n$
local affine coordinates. Set
\[
\check X(B):=\T^*_B/\check\Lambda.
\]
Of course $\T^*_B$ carries a canonical symplectic structure, and this
symplectic structure descends to $\check X(B)$.
\end{enumerate}
\qed
\end{definition}

We write $f:X(B)\rightarrow B$ and $\check f:\check X(B)\rightarrow B$
for these torus fibrations; these are clearly dual.

Now suppose in addition we have a Hessian metric $g$ on $B$,
with local potential function $K$. Then in fact
both $X(B)$ and $\check X(B)$ become K\"ahler manifolds:

\begin{proposition} $K\circ f$ is a (local) K\"ahler potential on
$X(B)$, defining a K\"ahler form $\omega=2i\partial\bar\partial(K\circ
f)$. This metric is
Ricci-flat if and only if $K$ satisfies
the real Monge-Amp\`ere equation
\[
\det {\partial^2 K\over \partial y_i\partial y_j}=constant.
\]
\end{proposition}

\begin{proof}
Working locally with affine coordinates $(y_i)$ and 
complex coordinates $z_j={1\over 2\pi i}\log q_j=x_j+i y_j$,
we compute $\omega=2i\partial\bar\partial(K\circ f)={i\over 2}
\sum {\partial^2 K\over \partial y_j\partial y_k} dz_j\wedge
d\bar z_k$
which is clearly positive. Furthermore,
if $\Omega=dz_1\wedge\cdots\wedge dz_n$, then
$\omega^n$ is proportional to $\Omega\wedge\bar\Omega$ if and only if
$\det (\partial^2 K/\partial y_j\partial y_k)$ is constant.
\end{proof}

We write this K\"ahler manifold as $X(B,K)$.

Dually we have

\begin{proposition}
In local canonical coordinates $y_i,\check x_i$ on $\T^*_B$, the functions
$z_j=\check x_j+i\partial K/\partial y_j$ on $\T^*_B$ induce a well-defined
complex structure on $\check X(B)$, with respect to which the canonical
symplectic form $\omega$ is a K\"ahler form of a metric. Furthermore
this metric is Ricci-flat if and only if $K$ satisfies the real
Monge-Amp\`ere equation
\[
\det {\partial^2 K\over \partial y_j\partial y_k}=constant.
\]
\end{proposition}

\begin{proof}
It is easy to see that an affine linear change in the coordinates
$y_j$ (and hence an appropriate change in the coordinates $\check x_j$)
results in a linear change of the coordinates $z_j$, so they induce
a well-defined complex structure invariant under $\check x_j\mapsto \check x_j+1$, and hence
a complex structure on $\check X(B)$. Then one computes that
\[
\omega=\sum d\check x_j\wedge dy_j={i\over 2}\sum g^{jk} dz_j\wedge d\bar z_k
\]
where $g_{ij}=\partial^2 K/\partial y_j\partial y_k$. Then the metric
is Ricci-flat if and only if $\det(g^{jk})=constant$, if and only if
$\det(g_{jk})=constant$.
\end{proof}

As before, we call this K\"ahler manifold $\check X(B,K)$.

This motivates the definition

\begin{definition} An affine manifold with metric of Hessian form
is a \emph{Monge-Amp\`ere manifold} if the local potential function $K$
satisfies the Monge-Amp\`ere equation $\det(\partial^2K/\partial y_i\partial
y_j)=constant$.
\end{definition}

Monge-Amp\`ere manifolds were first studied by Cheng and Yau in \cite{ChengYau}.

\begin{xca}
\label{caniso}
Show that the identification of $\T_B$ and $\T^*_B$ given by a Hessian metric
induces a canonical isomorphism $X(B,K)\cong\check X(\check B,\check K)$
of K\"ahler manifolds, where $(\check B,\check K)$ is the Legendre transform
of $(B,K)$.
\end{xca}

Finally, we note that a $B$-field can be introduced into this picture.
To keep life relatively simple (so as to avoid having to pass to generalized
complex structures \cite{HitGen}, \cite{Gual}, \cite{Oren}), 
we view the $B$-field as an element ${\bf B}\in
H^1(B,\Lambda_{\RR}/\Lambda)$, where $\Lambda_{\RR}=\Lambda\otimes_{\ZZ}
\RR$. Noting that a section of $\Lambda_{\RR}/\Lambda$ over an open
set $U$ can be viewed as a section of $\T_U/\Lambda|_U$, such a section
acts on $\T_U/\Lambda|_U$ via translation, and this action is in fact
holomorphic with respect to the standard semi-flat complex structure.
Thus a \v Cech 1-cocycle $(U_{ij},\beta_{ij})$
representing ${\bf B}$ allows us to reglue
$X(B)$ via translations over the intersections $U_{ij}$. This gives a new
complex manifold $X(B,{\bf B})$. If in addition there is a multi-valued
potential function $K$ defining a metric, these translations preserve
the metric and yield a K\"ahler manifold $X(B,{\bf B},K)$.

Thus the full toy version of mirror symmetry is as follows.
The data consists
of an affine manifold $B$ with potential $K$ and $B$-fields ${\bf B}
\in H^1(B,\Lambda_{\RR}/\Lambda)$, $\check {\bf B}\in H^1(B,
\check\Lambda_{\RR}/\check\Lambda)$. Now it is not difficult to see, and
you will have seen this already if you've done Exercise \ref{caniso}, that
the local system $\check\Lambda$ defined using the affine structure on $B$
is the same as the local system $\Lambda$ defined using the affine
stucture on $\check B$. So we say
the pair
\[
(X(B,{\bf B},K),\check{\bf B})
\]
is mirror to
\[
(X(\check B,\check {\bf B},\check K),\bf B).
\]

This provides a reasonably fulfilling picture of mirror symmetry
in a simple context. Many more aspects of mirror symmetry can be
worked out in this semi-flat context, see \cite{Leung}. However,
ultimately this only sheds limited insight into the general case.
The only compact Calabi-Yau manifolds with semi-flat Ricci-flat metric
which arise in this way are complex tori (shown by Cheng and Yau
in \cite{ChengYau}).
To deal with more interesting cases, we need to allow singular fibres,
and hence, singularities in the affine structure of $B$.

\section{Affine manifolds with singularities}

To deal with singular fibres, we define

\begin{definition}
A \emph{(tropical, integral) affine manifold with singularities}
is a $(C^0)$ manifold $B$ with an open subset $B_0\subseteq B$ which carries
a (tropical, integral) affine structure, and such that $\Gamma:=B
\setminus B_0$ is a locally finite union of locally closed submanifolds
of codimension $\ge 2$. 
\end{definition}

By way of example, let's explain how the Batyrev construction 
gives rise to a wide class of such manifolds. This construction is taken from
\cite{GBB}, where a more combinatorially complicated version is given
for complete intersections; see \cite{HZ} and \cite{HZ3} for an alternative
construction.

Let $\Delta$ be a reflexive polytope in $M_{\RR}=M\otimes_{\ZZ}\RR$,
where $M=\ZZ^n$; let $N$ be the dual lattice, $\nabla\subseteq
N_{\RR}$ the dual polytope given by
\[
\nabla:=\{n\in N_{\RR}|\hbox{$\langle m,n\rangle\ge -1$ for all $m\in\Delta$}
\}.
\] 
We assume $0\in\Delta$ is the unique interior lattice point of
$\Delta$. Let 
$\check\Sigma$ be the normal fan to $\nabla$, consisting of
cones over the faces of $\Delta$. Suppose we are
given a star subdivision of $\Delta$, with all vertices
being integral points, inducing a subdivision
$\check\Sigma'$ of the fan $\check\Sigma$. 
In addition suppose that
\[
\check h:M_{\RR}\rightarrow\RR
\]
is an (upper) strictly convex piecewise linear function on the fan
$\check\Sigma'$. Also, let 
\[
\check\varphi:M_{\RR}\rightarrow\RR
\]
be the piecewise linear function representing the anti-canonical class
of the toric variety $\PP_{\nabla}$;
i.e.\ $\check\varphi$ takes the value $1$ on the primitive generator of each
one-dimensional cone of $\check\Sigma$. Finally, assume that $\check h$
is chosen so that $\check h'=\check h-\check\varphi$
is a (not necessarily strictly) convex function.

Define, for any convex piecewise linear function $\check g$ on the fan
$\check\Sigma'$, the Newton polytope of $\check g$,
\[
\nabla^{\check g}:=\{n\in N_{\RR}|\hbox{$\langle m,n\rangle\ge -\check g(m)$ for
all $m\in M_{\RR}$}\}.
\]
In particular, 
\[
\nabla^{\check h}=\nabla^{\check h'}+\nabla^{\check\varphi}=
\nabla^{\check h'}+\nabla,
\]
where $+$ denotes Minkowski sum. Our goal will be to put an affine
structure with singularities on $B:=\partial\nabla^{\check h}$.
Our first method of doing this requires no choices.
Let $\P$ be the set of proper faces of $\nabla^{\check h}$. 
Furthermore, let
$\Bar(\P)$ denote the first barycentric subdivision of
$\P$ and let $\Gamma\subseteq B$ be the union of all simplices of
$\Bar(\P)$ not containing a vertex of $\P$ (a zero-dimensional cell)
or intersecting the interior of a maximal cell of $\P$. If we then
set $B_0:=B\setminus\Gamma$, we can define an affine structure on $B_0$
as follows. 
$B_0$ has an open cover
\[
\{W_{\sigma}|\hbox{$\sigma\in\P$ maximal}\}\cup
\{W_v|\hbox{$v\in\P$ a vertex}\}
\]
where $W_{\sigma}=\Int(\sigma)$, the interior of $\sigma$, and 
\[
W_v=\bigcup_{\tau\in\Bar(\P)\atop v\in\tau}\Int(\tau)
\]
is the (open) star of $v$ in $\Bar(\P)$. We define an affine chart
\[
\psi_{\sigma}:W_{\sigma}\rightarrow\AA^{n-1}\subseteq N_{\RR}
\]
given by the inclusion of $W_{\sigma}$ in
$\AA^{n-1}$, the affine hyperplane containing $\sigma$. Also, take
\[
\psi_v:W_v\rightarrow N_{\RR}/\RR v'
\]
to be the projection, where
$v$, being a vertex of $\nabla^{\check h}$, 
can be written uniquely as
$v'+v''$ with $v'$ a vertex of $\nabla$ and 
$v''$ a vertex of $\nabla^{\check h'}$.
One checks easily that for
$v\in\sigma$, $\psi_{\sigma}\circ\psi_v^{-1}$ is affine linear with
integral linear part (integrality follows from reflexivity of $\Delta$!)
so $B$ is a tropical affine manifold with singularities. 
Furthermore, if $\check h$ was chosen to have integral slopes, then
$B$ is integral.

We often would like to refine this construction, to get a finer polyhedral
decomposition $\P$ of $B$ and with it a somewhat more interesting 
discriminant locus $\Gamma$. One reason for doing so is that this
construction is clearly not mirror symmetric, as it depends only
on a star subdivision of $\Delta$ and not of $\nabla$. Furthermore,
a maximal star subdivision of $\nabla$ corresponds to what Batyrev
terms a MPCP (maximal projective crepant partial) resolution of $\PP_{\Delta}$,
and normally, we will wish to study hypersurfaces in a MPCP resolution
of $\PP_{\Delta}$ rather than in $\PP_{\Delta}$ itself.
To introduce this extra degree of flexibility, we need to make some choices,
which is done as follows.

First, choose
a star subdivision of $\nabla$, with all vertices being integral
points, inducing a refinement $\Sigma'$ of
the fan $\Sigma$ which is the normal fan to $\Delta$. 
This induces a polyhedral subdivision of
$\partial\nabla$, and we write the collection of cells
of this subdivision as $\P_{\partial\nabla}$. Note that because
$0\in \nabla$, we have
\[
\nabla^{\check h'}\subseteq \nabla^{\check h'}+\nabla=\nabla^{\check h}.
\]

\begin{definition} A subdivision $\P$ of
$\partial\nabla^{\check h}$ is \emph{good} with respect to 
$\P_{\partial\nabla}$ if it is induced by a subdivision 
$\P_{\nabla^{\check h}}$ of $\nabla^{\check h}$ satisfying the following
three properties:
\begin{enumerate}
\item $\nabla^{\check h'}$ is a union of cells in $\P_{\nabla^{\check h}}$.
\item All vertices of $\P_{\nabla^{\check h}}$ are contained either
in $\partial\nabla^{\check h}$ or in $\nabla^{\check h'}$.
\item Every cell $\sigma\in\P_{\nabla^{\check h}}$ with 
$\sigma\cap\partial\nabla^{\check h}\not=\emptyset$ and $\tau:=\sigma\cap
\partial\nabla^{\check h'}\not=\emptyset$ can be written as
\[
\sigma=(C(\sigma')+\tau)\cap \nabla^{\check h},
\]
with $\sigma'\in\P_{\partial\nabla}$ and $C(\sigma')$ the corresponding
cone in $\Sigma'$.
\end{enumerate}

If $\check h$ has integral slopes and all vertices of $\P_{\nabla^{\check h}}$
are integral, then we say $\P$ is integral.
\end{definition}
\vfill
\eject
The following picture shows what such a good subdivision may look like,
in the case that $\nabla$ is the Newton polytope of $\O_{\PP^2}(3)$:
\begin{center}
\includegraphics{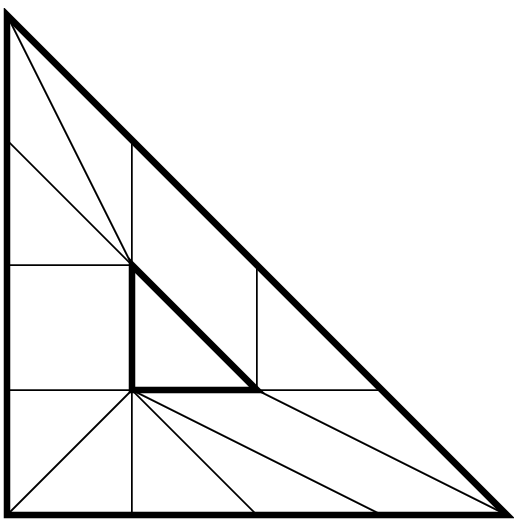}
\end{center}

Given a good decomposition $\P$ of 
\[
B:=\partial\nabla^{\check h},
\]
we once again obtain an affine structure with singularities on $B$, much
as before, defining the discriminant locus $\Gamma\subseteq B$
in terms of the first barycentric subdivision
$\Bar(\P)$ of this new polyhedral decomposition $\P$.
Then as before
$B_0:=B\setminus\Gamma$ has an open cover
\[
\{W_{\sigma}|\hbox{$\sigma\in\P$ maximal}\}\cup
\{W_v|\hbox{$v\in\P$ a vertex}\}
\]
where $W_{\sigma}=\Int(\sigma)$, the interior of $\sigma$, and 
\[
W_v=\bigcup_{\tau\in\Bar(\P)\atop v\in\tau}\Int(\tau)
\]
is the (open) star of $v$ in $\Bar(\P)$. We define an affine chart
\[
\psi_{\sigma}:W_{\sigma}\rightarrow\AA^{n-1}\subseteq N_{\RR}
\]
given by the inclusion of $W_{\sigma}$ in
$\AA^{n-1}$, the affine hyperplane containing $\sigma$. Also, take
$\psi_v:W_v\rightarrow N_{\RR}/\RR v'$ to be the projection, where
$v$ can be written uniquely as
$v'+v''$ with $v'$ an integral point of $\nabla$ and 
$v''\in\nabla^{\check h'}$.
As before, one checks easily that for
$v\in\sigma$, $\psi_{\sigma}\circ\psi_v^{-1}$ is affine linear with
integral linear part 
so $B$ is a tropical affine manifold with singularities. 
Furthermore, if $\check h$ has integral
slopes, and $\P$ is integral,
then the affine structure on $B$ is in fact
integral.

\begin{example}
\label{quintic}
Let $\Delta\subseteq\RR^4$ be the convex hull of the points
\[
(-1,-1,-1,-1), (1,0,0,0), (0,1,0,0), (0,0,1,0), (0,0,0,1),
\]
so $\nabla$ is the convex
hull of the points
\[
(-1,-1,-1,-1), (4,-1,-1,-1), (-1,4,-1,-1), (-1,-1,4,-1),
(-1,-1,-1,4).
\]
Take $\check h=\check\varphi$ and choose a star triangulation of 
$\nabla$. In this case $B=\partial\nabla$. It is
easy to see the affine structure on $B_0$ in fact extends across
the interior of all three-dimensional faces of $\nabla$. This
gives a smaller discriminant locus $\Gamma$ which, given a nice
regular triangulation of $\nabla$, looks like the following picture
in a neighbourhood of a $2$-face of $\nabla$: the light lines giving the 
triangulation and the dark lines the discriminant locus $\Gamma$.
\begin{center}
\includegraphics{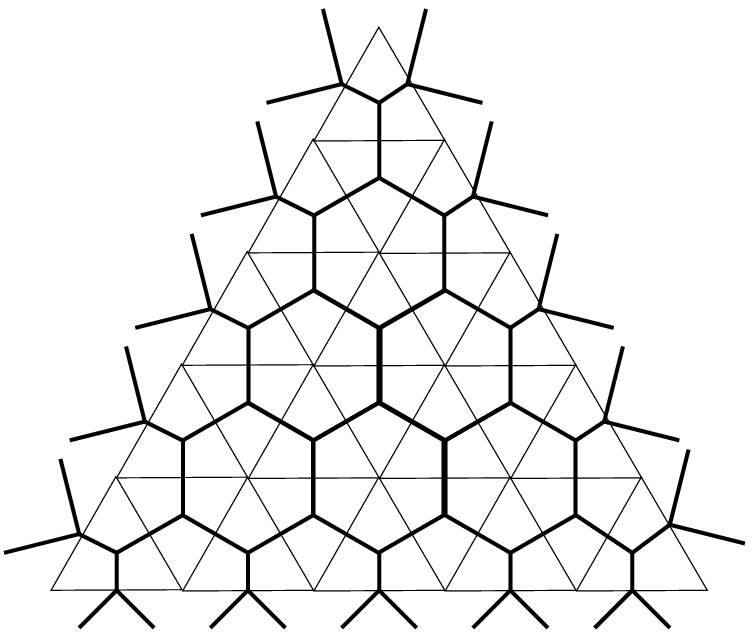}
\end{center}

In this picture, $\Gamma$ is a trivalent graph, with two types of trivalent
vertices. The ones along the edge are non-planar: the two additional
legs of $\Gamma$ drawn in this picture are contained in other two-faces
of $\nabla$. However $\Gamma$ is planar in the interior of this two-face.
\end{example}

In general, if the subdivisions $\Sigma'$ and $\check\Sigma'$ of $\Sigma$ and
$\check\Sigma$ respectively represent maximal projective crepant
partial resolutions of $\PP_{\Delta}$ and $\PP_{\nabla}$, and $\dim B=3$,
then only these sorts of trivalent vertices occur. More specifically, one
knows what the monodromy of the local systems $\Lambda$ and $\check\Lambda$
on $B_0$ look like at these vertices. If $v\in\Gamma$ is a vertex contained
in the interior of a two-face, then it is clear that the tangent space
to that two-face is invariant under parallel transport, in a neighbourhood
of $v$, of $\Lambda$.
A more careful analysis yields that the monodromy matrices for $\Lambda$ 
take the form, in a suitable basis,
\[
T_1=\begin{pmatrix} 1&0&0\\1&1&0\\0&0&1\end{pmatrix},
T_2=\begin{pmatrix} 1&0&0\\0&1&0\\1&0&1\end{pmatrix},
T_3=\begin{pmatrix} 1&0&0\\-1&1&0\\-1&0&1\end{pmatrix}.
\]
Here $T_1,T_2,T_3$ are given by parallel transport about loops around
the three edges of $\Gamma$ coming out of $v$. Of course, the monodromy
of $\check\Lambda$ is the transpose inverse of these matrices.
Similarly, if $v$ is a vertex of $\Gamma$ contained in an edge of 
$\nabla^{\check h}$, then the monodromy will take the form
\[
T_1=\begin{pmatrix} 1&-1&0\\0&1&0\\0&0&1\end{pmatrix},
T_2=\begin{pmatrix} 1&0&-1\\0&1&0\\0&0&1\end{pmatrix},
T_3=\begin{pmatrix} 1&1&1\\0&1&0\\0&0&1\end{pmatrix}.
\]
So we see that the monodromy of the two types of vertices are interchanged
by looking at $\Lambda$ and $\check\Lambda$. 

One main result of \cite{TMS} is

\begin{theorem}
If $B$ is a three-dimensional tropical
affine manifold with singularities such that $\Gamma$ is
trivalent and the monodromy of $\Lambda$ at each vertex is one of the
above two types, then $f_0:X(B_0)\rightarrow B_0$ can be compactified to
a topological fibration $f:X(B)\rightarrow B$. Dually, $\check f_0:\check X(B_0)
\rightarrow B_0$ can be compactified to a topological fibration
$\check f:\check X(B)\rightarrow B$.
\end{theorem}

We won't give any details here of how this is carried out, but it is
not particularly difficult, as long as one restricts to the category
of topological (not $C^{\infty}$) manifolds. However, it is interesting
to look at the singular fibres we need to add in this compactification.

If $b\in\Gamma$ is a point which is not a vertex of $\Gamma$, then $f^{-1}(b)$
is homeomorphic to $I_1\times S^1$, where $I_1$ denotes a Kodaira type $I_1$
elliptic curve, i.e.\ a pinched torus.

If $b$ is a vertex of $\Gamma$, with monodromy of the first type, then
$f^{-1}(b)=S^1\times S^1\times S^1/\sim$, with $(a,b,c)\sim (a',b',c')$
if $(a,b,c)=(a',b',c')$ or $a=a'=1$, where $S^1$ is identified with the unit
circle in $\CC$.
This is the three-dimensional analogue
of a pinched torus, and $\chi(f^{-1}(b))=+1$. We call this a \emph{positive}
fibre.

If $b$ is a vertex of $\Gamma$, with monodromy of the second type, then
$f^{-1}(b)$ can be described as $S^1\times S^1\times S^1/\sim$,
with $(a,b,c)\sim (a',b',c')$ if $(a,b,c)=(a',b',c')$ or $a=a'=1$, $b=b'$,
or $a=a',b=b'=1$.
The singular locus of this fibre is a figure eight, and $\chi(f^{-1}(b))=-1$.
We call this a \emph{negative} fibre.

So we see a very concrete local consequence of SYZ duality: namely
in the compactifications $X(B)$ and $\check X(B)$, the positive and
negative fibres are interchanged. Of course, this results in the
observation that Euler characteristic changes sign under mirror symmetry
for Calabi-Yau threefolds.

\begin{example}
Continuing with Example \ref{quintic}, it was proved in \cite{TMS} 
that $\check X(B)$ is homeomorphic to the quintic and $X(B)$
is homeomorphic to the mirror quintic.
\end{example}

\begin{remark}
Haase and Zharkov in \cite{HZ} gave a different description of what is
the same affine structure. Their construction has the advantage
that it is manifestly dual. In other words, in our construction,
we can interchange the role of $\Delta$ and $\nabla$ to get two different
affine manifolds, with $B_{\Delta}$
the affine manifold with singularities structure on $\partial\Delta^{h}$
and $B_{\nabla}$ the affine manifold with singularities structure
on $\partial\nabla^{\check h}$. It is not obvious that these are ``dual''
affine manifolds, at least in the sense that $X(B_{\nabla})$ is homeomorphic
to $\check X(B_{\Delta})$ and $\check X(B_{\nabla})$ is homeomorphic
to $X(B_{\Delta})$. In the construction given above, this follows from
the discrete Legendre transform we will discuss in \S 7. On the other hand,
the construction I give here will arise naturally from the degeneration
construction discussed later in this paper.

Ruan in \cite{Ruan} gave a description of \emph{Lagrangian} torus fibrations
for hypersurfaces in toric varieties using a symplectic flow argument, 
and his construction should
coincide with a \emph{symplectic} compactification of the symplectic
manifolds $\check X(B_0)$. In the three-dimensional case, such a 
symplectic compactification has now been constructed by Ricardo
Casta\~no-Bernard and Diego Matessi \cite{CastMat}. If this compactification
is applied to the affine manifolds with singularities described here,
the resulting symplectic manifolds should be symplectomorphic
to the corresponding toric hypersurface, but this has not yet been shown.

I should also point out that the explicit compactifications mentioned in
three dimensions can be carried out in all dimensions, and will be done
so in \cite{tori}. 
We will show there in a much more general context that
these compactifications are then homeomorphic to the expected Calabi-Yau
manifolds.
\end{remark}

\section{The problems with the SYZ conjecture, and how to get around them}

The previous section demonstrates that the SYZ conjecture gives a beautiful
description of mirror symmetry at a purely topological level. This, by itself,
can often be useful, but unfortunately is not strong enough to get
at really interesting aspects of mirror symmetry, such as instanton
corrections. For a while, though, many of us were hoping that the strong
version of duality we have just seen would hold at the special Lagrangian
level. This would mean that a mirror pair $X,\check X$ would
possess special Lagrangian torus fibrations
$f:X\rightarrow B$ and $\check f:\check X\rightarrow B$ with codimension
two discriminant locus, and the discriminant loci of $f$ and
$\check f$ would coincide.
These fibrations would then be dual away from the discriminant locus.

There are examples of special Lagrangian fibrations on non-compact
toric varieties $X$ with this behaviour. In particular, if $\dim X=n$
with a $T^{n-1}$ action on $X$ preserving the holomorphic $n$-form,
and if $X$ in addition carries a Ricci-flat metric which is invariant under
this action, then $X$ will have a very nice special Lagrangian fibration
with codimension two discriminant locus. (See \cite{SLAGex} and 
\cite{Gold}). However, Dominic Joyce (\cite{Joyce} and other
papers cited therein) began studying some 
three-dimensional $S^1$-invariant
examples, and discovered quite different behaviour. There is an argument
that if a special Lagrangian fibration is $C^{\infty}$, then the
discriminant locus will be (Hausdorff) codimension two. However, Joyce
discovered examples which were not differentiable, but only piecewise
differentiable, and furthermore, had a codimension one discriminant locus:

\begin{example}
Define $F:\CC^3\rightarrow \RR\times\CC$ by
$F(z_1,z_2,z_3)=(a,c)$ with
$2a=|z_1|^2-|z_2|^2$ and
\[
c=\begin{cases}
z_3&a=z_1=z_2=0\\
z_3-\bar z_1\bar z_2/|z_1|& a\ge 0, z_1\not=0\\
z_3-\bar z_1\bar z_2/|z_2|&a<0.
\end{cases}
\]
It is easy to see that if $a\not=0$, then $F^{-1}(a,c)$ is homeomorphic to
$\RR^2\times S^1$, while if $a=0$, then $F^{-1}(a,c)$ is a cone over $T^2$:
essentially, one copy of $S^1$ in $\RR^2\times S^1$ collapses to a point. 
In addition, all fibres of this map are special Lagrangian, and it is obviously
only piecewise smooth. The discriminant locus is the entire plane given
by $a=0$.
\end{example}

This example forces a reevaluation of the strong form of the SYZ conjecture.
In further work Joyce found evidence for a more likely picture
for general special Lagrangian fibrations in three dimensions. The discriminant
locus, instead of being a codimension two graph, will be a codimension one
blob. Typically the union of the singular points of singular fibres will
be a Riemann surface, and it will map to an amoeba shaped set in $B$, i.e.\
the discriminant locus looks like the picture on the right rather than the left,
and will be a fattening of the old picture of a codimension two discriminant.

\begin{center}
\includegraphics{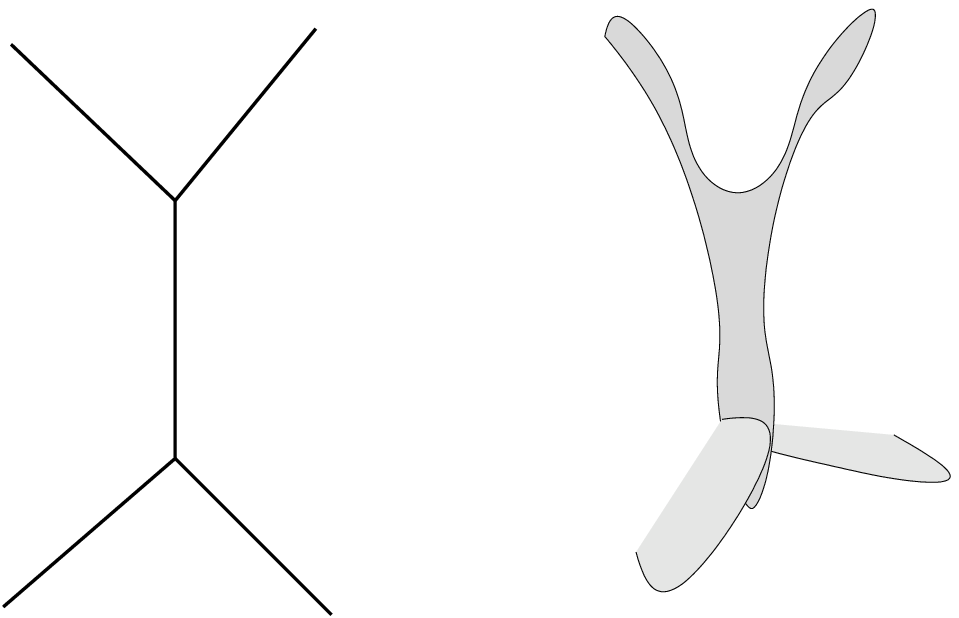}
\end{center}

Joyce made some additional arguments to suggest that this fattened 
discriminant locus must look fundamentally different in a neighbourhood
of the two basic types of vertices we saw in the previous section,
with the two types of vertices expected to appear pretty much as depicted in
the above picture.
Thus the strong form of duality mentioned above, where we expect the
discriminant loci of the special Lagrangian fibrations on a mirror pair
to be the same, cannot hold. If this is the case, one needs to replace
this strong form of duality with a weaker form.

It seems likely that the best way to rephrase the SYZ conjecture is
in a limiting form. Mirror symmetry as we currently understand it has to
do with degenerations of Calabi-Yau manifolds. Given a flat family $f:\X
\rightarrow D$ over a disk $D$, with the fibre $\X_0$ over $0$
singular and all other fibres $n$-dimensional
Calabi-Yau manifolds, we say the family is \emph{maximally unipotent}
if the monodromy transformation $T:H^n(\X_t,\QQ)\rightarrow H^n(\X_t,\QQ)$
($t\in D$ non-zero) satisfies $(T-I)^{n+1}=0$ but $(T-I)^n\not=0$.
It is a standard fact of mirror symmetry that mirrors should be associated
to maximally unipotent degenerations of Calabi-Yau manifolds. In particular,
given two different maximally unipotent degenerations in a single complex
moduli space for some Calabi-Yau manifold, one might obtain different
mirror manifolds. Sometimes these different mirror manifolds are birationally
equivalent, as studied in \cite{AGM}, or are genuinely different,
see \cite{Rodland}. 

We recall the definition of Gromov-Hausdorff convergence,
a notion of convergence of a sequence of metric spaces.

\begin{definition}
Let $(X,d_X)$, $(Y,d_Y)$ be two compact
metric spaces. Suppose there exists maps $f:X\rightarrow Y$
and $g:Y\rightarrow X$ (not necessarily continuous) such that
for all $x_1,x_2\in X$,
$$|d_X(x_1,x_2)-d_Y(f(x_1),f(x_2))|<\epsilon$$
and for all $x\in X$,
$$d_X(x,g\circ f(x))<\epsilon,$$
and the two symmetric properties for $Y$ hold. Then we say the
Gromov--Hausdorff distance between $X$ and $Y$ is at most $\epsilon$.
The Gromov--Hausdorff distance $d_{GH}(X,Y)$ is the infimum of all
such $\epsilon$.
\end{definition}

It follows from results of Gromov (see for example \cite{Petersen}, 
pg. 281, Cor.
1.11) that the space of compact Ricci-flat manifolds with diameter $\le C$
is precompact with respect to Gromov-Hausdorff distance, 
i.e.\ any sequence of such manifolds has a subsequence
converging with respect to the Gromov-Hausdorff distance to a metric
space. This metric space could be quite bad; this is quite outside
the realm of algebraic geometry! Nevertheless, this raises the following
natural question. Given a maximally unipotent
degeneration of Calabi-Yau manifolds $\X\rightarrow D$, take a 
sequence $t_i\in D$ converging to $0$, and consider a sequence $(\X_{t_i},
g_{t_i})$, where $g_{t_i}$ is a choice of Ricci-flat metric chosen so
that $Diam(g_{t_i})$ remains bounded. What is the Gromov-Hausdorff
limit of $(\X_{t_i},g_{t_i})$, or the limit of some convergent subsequence?

\begin{example} Consider a degenerating family of elliptic curves, with
periods $1$ and ${1\over 2\pi i}\log t$. If we take $t$ approaching
$0$ along the positive real axis, then we can just view this as a family
of elliptic curves $\X_{\alpha}$
with period $1$ and $i\alpha$ with $\alpha\rightarrow\infty$.
If we take the standard Euclidean metric $g$ on $\X_{\alpha}$, then
the diameter of $\X_{\alpha}$ is unbounded. To obtain a bounded diameter,
we replace $g$ by $g/\alpha^2$; equivalently, we can keep $g$ fixed
on $\CC$ but change the periods of the elliptic curve to $1/\alpha, i$.
It then becomes clear that the Gromov-Hausdorff limit of such a sequence
of elliptic curves is a circle $\RR/\ZZ$.
\end{example}

This simple example motivates the first conjecture about maximally
unipotent degenerations, conjectured independently by myself and Wilson
on the one hand \cite{GrWi} and Kontsevich and Soibelman \cite{KS} on the other.

\begin{conjecture}
\label{GWKSconj}
Let $\X\rightarrow D$ be a maximally unipotent
degeneration of simply-connected 
Calabi-Yau manifolds with full $SU(n)$ holonomy, 
$t_i\in D$ with $t_i\rightarrow 0$,
and let $g_i$ be a Ricci-flat metric on $\X_{t_i}$ normalized to have
fixed diameter $C$. Then a convergent subsequence of $(\X_{t_i},g_i)$
converges to a metric space $(X_{\infty},d_{\infty})$, where $X_{\infty}$
is homeomorphic to $S^n$. Furthermore, $d_{\infty}$ is induced by
a Riemannian metric on $X_{\infty}\setminus\Gamma$, where $\Gamma\subseteq
X_{\infty}$ is a set of codimension two.
\end{conjecture}

Here the topology of the limit depends on the nature of the non-singular
fibres $\X_t$; for example,
if instead $\X_t$ was hyperk\"ahler, then we would expect
the limit to be a projective space. Also, even in the case of full
$SU(n)$ holonomy, if $\X_t$ is not simply connected, we would expect
limits such as $\QQ$-homology spheres to arise.

Conjecture \ref{GWKSconj} is directly inspired by the SYZ conjecture. Suppose
we had special Lagrangian fibrations $f_i:\X_{t_i}\rightarrow B_i$.
Then as the maximally unipotent degeneration is approached, it is possible
to see that the volume of the fibres of these fibrations
go to zero. This would suggest these fibres collapse, hopefully
leaving the base as the limit.

This conjecture was proved by myself and Wilson
for K3 surfaces in \cite{GrWi}. The proof relies
on a detailed analysis of the behaviour of Ricci-flat metrics in the limit,
and also on the existence of explicit local models for Ricci-flat metrics
near singular fibres of special Lagrangian fibrations. 

The motivation for this conjecture from SYZ also provides a limiting form
of the conjecture. There are any number of problems with trying to prove
the existence of special Lagrangian fibrations on Calabi-Yau manifolds.
Even the existence of a single special Lagrangian torus near a maximally
unipotent degeneration is unknown, but we expect it should be easier to
find them as we approach the maximally unipotent point. Furthermore,
even if we find a special Lagrangian torus, we know that it moves
in an $n$-dimensional family, but we don't know its deformations fill
out the entire manifold. In addition, there is no guarantee that even if it
does, we obtain a foliation of the manifold: nearby special Lagrangian
submanifolds may intersect. (For an example, see \cite{Matessi}.)
So instead, we will just look at the moduli space
of special Lagrangian tori.

Suppose, given $t_i\rightarrow 0$, that for $t_i$ sufficiently close
to zero, there is a special Lagrangian $T^n$ whose homology class is
invariant under monodromy, or more specifically, generates the space $W_0$
of the monodromy weight filtration
(this is where we expect to find fibres
of a special Lagrangian fibration associated to a maximally unipotent
degeneration). Let $B_{0,i}$ be the moduli space of deformations
of this torus; every point of $B_{0,i}$ corresponds to a smooth special
Lagrangian torus in $\X_{t_i}$. This manifold then comes
equipped with the McLean metric and affine structures
defined in \S 2. One can then
compactify $B_{0,i}\subseteq B_i$, (probably by taking the closure
of $B_{0,i}$ in the space of special Lagrangian currents; the details
aren't important here). This gives a series of metric spaces $(B_i,d_i)$
with the metric $d_i$ induced by the McLean metric. If the McLean
metric is normalized to keep the diameter of $B_i$ constant independent of $i$,
then we can hope that $(B_i,d_i)$ converges to a compact metric space
$(B_{\infty},d_{\infty})$. Here then is the limiting form of SYZ:

\begin{conjecture} If $(\X_{t_i},g_i)$ converges to $(X_{\infty},g_{\infty})$
and $(B_i,d_i)$ is non-empty for large $i$ and converges to
$(B_{\infty},d_{\infty})$, then $B_{\infty}$ and $X_{\infty}$
are isometric up to scaling. Furthermore, there is a subspace $B_{\infty,0}
\subseteq B_{\infty}$ with $\Gamma:=B_{\infty}\setminus B_{\infty,0}$ 
of Hausdorff codimension 2 in $B_{\infty}$
such that $B_{\infty,0}$ is a Monge-Amp\`ere manifold, with the Monge-Amp\`ere
metric inducing $d_{\infty}$ on $B_{\infty,0}$.
\end{conjecture}

Essentially what this is saying is that as we approach the maximally
unipotent degeneration, we expect to have a special Lagrangian fibration
on larger and larger subsets of $\X_{t_i}$. Furthermore, in the limit,
the codimension one discriminant locus suggested by Joyce converges to
a codimension two discriminant locus, and (the not necessarily
Monge-Amp\`ere, see \cite{Matessi}) Hessian metrics on $B_{0,i}$
converge to a Monge-Amp\`ere metric.

The main point I want to get at here is that
it is likely the SYZ conjecture is only ``approximately'' correct, and
one needs to look at the limit to have a hope of proving anything.
On the other hand, the above conjecture seems likely to be accessible
by currently understood techniques, though with a lot of additional work,
and I wouldn't be surprised to see it proved in the next few years.

How do we do mirror symmetry using this modified version of the SYZ conjecture?
Essentially, we would follow these steps:
\begin{enumerate}
\item 
We begin with a maximally unipotent degeneration of Calabi-Yau manifolds
$\X\rightarrow D$, along with a choice of polarization. This gives
us a K\"ahler class $[\omega_t]\in H^2(\X_t,\RR)$ for each $t\in D\setminus
0$, represented by $\omega_t$ the K\"ahler form of a Ricci-flat metric
$g_t$. 
\item Identify the Gromov-Hausdorff limit of a sequence $(\X_{t_i}, r_ig_{t_i})$
with $t_i\rightarrow 0$, and $r_i$ a scale factor which keeps the diameter
of $\X_{t_i}$ constant. The limit will be, if the above conjectures work,
an affine manifold with singularities $B$ along with a Monge-Amp\`ere metric.
\item Perform a Legendre transform to obtain a new affine manifold with
singularities $\check B$, though with the same metric.
\item Try to construct a compactification of $X_{\epsilon}(\check B_0)$
for small $\epsilon>0$ to obtain a complex manifold $X_{\epsilon}(\check B)$.
This will be the mirror manifold.
\end{enumerate}

Actually, we need to elaborate on this last step a bit more. The
problem is that while we expect that it should be possible in general
to construct symplectic compactifications of the symplectic manifold
$\check X(B_0)$ (and hence get the mirror as a symplectic manifold), 
we don't expect to be able to compactify $X_{\epsilon}(\check B_0)$
as a complex manifold. Instead, the expectation is that a small deformation
of $X_{\epsilon}(\check B_0)$ is necessary before it can be compactified.
Furthermore, this small deformation is critically important in mirror
symmetry: \emph{it is this small deformation which provides the $B$-model
instanton corrections}. 

Because of the importance of this last issue, it has already been studied
by several authors: Fukaya in \cite{Fukaya} has studied the problem directly using
heuristic ideas, while Kontsevich and Soibelmann \cite{KS2} have modified the
problem of passing from an affine manifold to a complex manifold by instead
producing a non-Archimedean space. We will return to these issues later
in this paper, when I discuss my own work with Siebert which has been
partly motivated by the same problem. Because this last item is so 
important, let's give it a name:

\begin{question}[The reconstruction problem, Version I]
\label{reconstruct1}
Given a tropical affine manifold with singularities $B$, construct
a complex manifold $X_{\epsilon}(B)$ which is a compactification
of a small deformation of $X_{\epsilon}(B_0)$.
\end{question} 

I do not wish to dwell further on this version of the SYZ
conjecture here, because it lies mostly
in the realm of analysis and differential geometry and the behaviour of
Ricci-flat metrics, and will give us little insight into what makes
subtler aspects of traditional mirror symmetry work: for example, how exactly
do instanton corrections arise? So we move on to explain how the
limiting form of SYZ inspired a more algebro-geometric form of SYZ,
which in turn avoids all analytic problems and holds out great promise
for understanding the fundamental mysteries of mirror symmetry.

\section{Gromov-Hausdorff limits, algebraic degenerations, and mirror
symmetry}

We now have two notions of limit: the familiar algebro-geometric 
notion of a flat degenerating family $\X\rightarrow D$ over a disk
on the one hand, and the Gromov-Hausdorff limit on the other. Kontsevich
had an important insight (see \cite{KS}) into the connection between these
two. In this section I will give a rough idea of how and why this works.

Very roughly speaking, the Gromov-Hausdorff limit $(\X_{t_i},g_{t_i})$
as $t_i\rightarrow 0$ should coincide, topologically, with the dual
intersection complex of the singular fibre $\X_0$. More precisely,
in a relatively simple situation,
suppose $f:\X\rightarrow D$ is relatively minimal (in the sense of Mori)
and normal crossings,
with $\X_0$ having irreducible components $X_1,\ldots,X_m$. The dual 
intersection complex of $\X_0$ is the simplicial complex with vertices
$v_1,\ldots,v_m$, and which contains a simplex $\langle v_{i_0},\ldots,
v_{i_p}\rangle$ if $X_{i_0}\cap\cdots\cap X_{i_p}\not=\emptyset$.

Let us explain roughly why this should be, first by looking at a standard
family of degenerating elliptic curves with periods $1$ and ${n\over 2\pi i}
\log t$ for $n$ a positive integer. Such a family over the punctured disk
is extended to a family over the disk by adding an $I_n$ (a cycle of $n$
rational curves) fibre over the origin.

Taking a sequence $t_i\rightarrow 0$ with $t_i$ real and positive gives
a sequence of elliptic curves of the form $X_{\epsilon_i}(B)$ where:
$B=\RR/n\ZZ$ and
$\epsilon_i=-{2\pi\over\ln t_i}$. In addition, the metric on 
$X_{\epsilon_i}(B)$, properly scaled, comes from the constant Hessian
metric on $B$. So we wish to explain how $B$ is related to the geometry 
near the singular fibre. To this end,
let $X_1,\ldots,X_n$ be the irreducible components of $\X_0$; these
are all $\PP^1$'s. Let $P_1,\ldots,P_n$ be the singular points of $\X_0$.

We'll consider two sorts of open sets in $\X$. For the first type,
choose a coordinate $z$ on $X_i$, with $P_i$ given by $z=0$ and
$P_{i+1}$ given by $z=\infty$. Let $U_i\subseteq D_i$ be the open set
$\{z|\delta\le |z| \le 1/\delta\}$ for some small fixed $\delta$. Then one
can find a neighbourhood $\tilde U_i$ of $U_i$ in $\X$ such that 
$\tilde U_i$ is biholomorphic to $U_i\times D_{\rho}$ for $\rho>0$ sufficiently
small, $D_{\rho}$ a disk of radius $\rho$ in $\CC$, and
$f|_{\tilde U_i}$ is the projection onto $D_{\rho}$. 

On the
other hand, each $P_i$ has a neighbourhood $\tilde V_i$ in $\X$
biholomorphic to a polydisk
$\{(z_1,z_2)\in\CC^2||z_1|\le \delta', |z_2|\le\delta'\}$ on which $f$
takes the form $z_1z_2$. 

If $\delta$ and $\delta'$ are chosen correctly,
then for $t$ sufficiently close to zero,
\[
\{\tilde V_i\cap\X_t|1\le i\le n\}\cup \{\tilde U_i\cap\X_t|1\le i\le n\}
\]
form an open cover of $\X_t$. Now each of the sets in this open cover
can be written as $X_{\epsilon}(U)$ for some $U$ a one-dimensional
(non-compact) affine manifold and
$\epsilon=-2\pi/\ln|t|$. If $U$ is an open interval $(a,b)\subseteq
\RR$, then $X_{\epsilon}(U)$ is biholomorphic to the annulus
\[
\{z\in\CC| e^{-2\pi b/\epsilon}\le |z|\le e^{-2\pi a/\epsilon}\}
\]
as $q=e^{2\pi i(x+i y)/\epsilon}$ is a holomorphic coordinate on
$X_{\epsilon}((a,b))$. 
Thus 
\[
\tilde U_i\cap \X_t\cong X_{\epsilon}\left(\left({\epsilon\ln\delta\over 2\pi},
-{\epsilon\ln\delta\over 2\pi}\right)\right)
\]
with $\epsilon=-2\pi/\ln|t|$. As
$t\rightarrow 0$, the interval $(\epsilon\ln\delta/2\pi,
-\epsilon\ln\delta/2\pi)$ shrinks to a point. So $\tilde U_i\cap
\X_t$ is a smaller and smaller open subset of $\X_t$ as $t\rightarrow 0$
when we view things in this way. This argument suggests that every
irreducible component should be associated to a point on $B$.

Now look at $\tilde V_i\cap\X_t$. This is
\begin{eqnarray*}
\{(z_1,z_2)\in\CC^2||z_1|,|z_2|<\delta', z_1z_2=t\}
&\cong&\{z\in\CC||t|/\delta'\le |z|\le \delta'\}\\
&\cong& X_{\epsilon}\left({-\epsilon\over 2\pi}\ln\delta',
{\epsilon\over 2\pi} (\ln\delta'-\ln |t|)\right)
\end{eqnarray*}
with $\epsilon=-2\pi/\ln|t|$. This interval approaches the unit interval
$(0,1)$ as $t\rightarrow 0$. So the open set $\tilde V_i\cap \X_t$ ends up being
a large portion of $\X_t$. We end up with $\X_t$, for small $t$, being
a union of open sets of the form 
$X_{\epsilon}((i+\epsilon',i+1-\epsilon'))$ (i.e.\ $\tilde V_i\cap\X_{\epsilon}$)
and $X_{\epsilon}((i-\epsilon'',i+\epsilon''))$ (i.e.\ $\tilde U_i\cap
\X_t$) for $\epsilon'$, $\epsilon''$ sufficiently small. These should glue,
at least approximately, to give $X_{\epsilon}(B)$. So we see
that irreducible components of $\X_0$ seem to coincide with points on $B$,
but intersections of components
coincide with lines. In this way we see the dual
intersection complex emerge.

\medskip

Let us make one more observation before beginning with rigorous
results in the next section. Suppose more generally we had a \emph{Gorenstein
toroidal
crossings} degeneration of Calabi-Yau manifolds $f:\X\rightarrow D$. 
This means that every point $x\in\X$ has a neighbourhood isomorphic
to an open set in an affine Gorenstein (i.e.\ the canonical class is
a Cartier divisor) toric variety, with $f$ given
locally by a monomial which vanishes exactly
to order $1$ on each codimension one
toric stratum. This is a generalization of the notion of normal crossings,
see \cite{ss}.
Very roughly, the above argument suggests that each irreducible component
of the central fibre will correspond to a point of the Gromov-Hausdorff
limit. The following exercise shows what kind of contribution to $B$
to expect from a point $x\in\X_0$ which is a zero-dimensional stratum in
$\X_0$.

\begin{xca}
\label{gorensteinlimit}
Suppose there is a point $x\in\X_0$ which has
a neighbourhood isomorphic to a neighbourhood of a dimension zero torus
orbit of an affine Gorenstein toric variety $Y_x$. 
Such an affine variety is
specified as follows.
Set $M=\ZZ^n$, $M_{\RR}=M\otimes_{\ZZ}\RR$, $N=\Hom_{\ZZ}(M,\ZZ)$,
$N_{\RR}=N\otimes_{\ZZ}\RR$ as in \S 4. Then
there is a lattice polytope $\sigma\subseteq M_{\RR}$, 
$n=\dim\X_t$, $C(\sigma):=\{(rm,r)| m\in\sigma,r\ge 0\}\subseteq
M_{\RR}\oplus\RR$, $P:=\dual{C(\sigma)}\cap (N\oplus\ZZ)$ the monoid
determined by the dual of the cone $C(\sigma)$, and finally, $Y_x
=\Spec \CC[P]$, and $f$ coincides with the monomial $z^{(0,1)}$.

Let us now take a small neighbourhood of $x$ of the form
\[
\tilde 
U_{\delta}=\{y\in \Spec \CC[P]\,|\,\hbox{$|z^p|<\delta$ for all $p\in P$}\}.
\]
This is an open set as the condition $|z^p|<\delta$ can be tested on a
finite generating set for $P$, provided that $\delta<1$. 
Then show that for a given $t$, $|t|<1$ and $\epsilon=-2\pi/\log|t|$, if 
\[
\sigma_t:=\{m\in M_{\RR}|\hbox{$\langle p,(m,1)\rangle>{\log\delta\over
\log |t|}$ for all $p\in P$}\},
\]
then 
\[
f^{-1}(t)\cap \tilde U_{\delta}\cong X_{\epsilon}(\sigma_t).
\]
Note that 
\[
\sigma:=\{m\in M_{\RR}|\hbox{$\langle p,(m,1)\rangle\ge 0$ for all $p\in P$}\},
\]
so $\sigma_t$ is an open subset of $\sigma$, and
as $t\rightarrow 0$, $\sigma_t$ converges to the interior of
$\sigma$. 
\qed
\end{xca}

This observation will hopefully motivate the basic construction
of the next section.

\section{Toric degenerations, the intersection complex and its dual}

We now return to rigorous statements. I would like to explain the basic
ideas behind the program launched in \cite{PartI}. While I will use
the previous sections as motivation, this work actually got its start
when Siebert began a program of studying mirror symmetry via degenerations
of Calabi-Yau manifolds. Work of Schr\"oer and Siebert \cite{ssKod}, \cite{ss}
led Siebert
to the idea that log structures on degenerations of Calabi-Yau manifolds
would allow one to view mirror symmetry as an operation performed on
degenerate Calabi-Yau varieties. Siebert observed that at a combinatorial
level, mirror symmetry exchanged data pertaining to the log structure
and a polarization. This will be explained more clearly in the following
section, where we introduce log structures. Together, Siebert and I realised
that the combinatorial data he was considering could be encoded naturally
in the dual intersection complex of the degeneration, and that mirror
symmetry then corresponded to a discrete Legendre transform on the dual
intersection complex. It then became apparent that this approach provided
an algebro-geometrization of the SYZ conjecture.

Here I will explain this program from the opposite direction, starting with
the motivation of the previous section for introducing the dual intersection
complex, and then work backwards until we arrive naturally at log structures.
Much of the material in this section comes from \cite{PartI}, \S 4. 

\begin{definition}
\label{toricdegen}
Let $f:\X\rightarrow D$ be a proper flat family of relative dimension
$n$, where $D$ is a disk and $\X$ is a complex analytic space
(not necessarily non-singular). We say $f$ is a {\it toric degeneration}
of Calabi-Yau varieties if
\begin{enumerate}
\item
$\X_t$ is an irreducible normal Calabi-Yau variety with only
canonical singularities for $t\not=0$. (The reader may
like to assume $\X_t$ is smooth for $t\not=0$).
\item
If $\nu:\tilde\X_0\to\X_0$ is the normalization,
then $\tilde\X_0$ is a disjoint union of toric varieties,
the conductor locus $C\subseteq\tilde\X_0$ is reduced,
and the map $C\to\nu(C)$ is unramified and generically
two-to-one. (The conductor locus is a naturally defined scheme
structure on the set where $\nu$ is not an isomorphism.) The square
\[\begin{CD}
C@>>> \tilde\X_0\\
@VVV @VV{\nu}V\\
\nu(C)@>>> \X_0
\end{CD}\]
is cartesian and cocartesian.
\item $\X_0$ is a reduced Gorenstein space
and the conductor locus $C$
restricted to each irreducible component of $\tilde\X_0$ is the union
of all toric Weil divisors of that component. 
\item There exists a closed subset $Z\subseteq\X$ of relative
codimension $\ge 2$ such that $Z$ satisfies the following properties:
$Z$ does not contain the image under $\nu$
of any toric stratum of $\tilde\X_0$,
and for any point $x\in \X\setminus Z$, there is a neighbourhood
$\tilde U_x$ (in the analytic topology) of $x$, an
$n+1$-dimensional affine toric variety $Y_x$, a regular function
$f_x$ on $Y_x$ given by a monomial, and a commutative diagram
$$\begin{matrix}
\tilde U_x&\mapright{\psi_x}&Y_x\cr
\mapdown{f|_{\tilde U_{x}}}&&\mapdown{f_x}\cr
D'&\mapright{\varphi_x}&\CC\cr
\end{matrix}$$
where $\psi_x$ and $\varphi_x$ are open embeddings and $D'\subseteq D$. 
Furthermore,
$f_x$ vanishes precisely once on each toric divisor of $Y_x$.
\end{enumerate}
\end{definition}

\begin{example}
Take $\X$ to be defined by the equation $tf_4+z_0z_1z_2z_3=0$
in $\PP^3\times D$, where $D$ is a disk with coordinate
$t$ and $f_4$ is a general homogeneous quartic polynomial on $\PP^3$.
It is easy to see that $\X$ is singular at the locus
\[
\{t=f_4=0\}\cap Sing(\X_0).
\]
As $\X_0$ is the coordinate tetrahedron, the singular locus of
$\X_0$ consists of the 
six coordinate lines of $\PP^3$, and $\X$ has four singular points along
each such line, for a total of 24 singular points. 
Take $Z=Sing(\X)$. Then away from
$Z$, the projection $\X\rightarrow D$ is normal crossings, which
yields condition (4) of the definition of toric degeneration. It is easy
to see all other conditions are satisfied.
\end{example}

\begin{example}
\label{toricdegenexample2}
Let $\Delta\subseteq M_{\RR}$ be a reflexive polytope with dual $\nabla
\subseteq N_{\RR}$. Choose an \emph{integral} strictly convex
piecewise linear function
$\check h:M_{\RR}\rightarrow\RR$
as in \S 4. Consider the family in $\PP_{\Delta}\times D$ defined by
\begin{equation}
\label{toricdegenexam}
z^0+t\sum_{m\in\Delta\cap M} a_m t^{\check h'(m)}z^m=0.
\end{equation}
Here $z^m$ denotes the section of $\O_{\PP_{\Delta}}(1)$ determined by
the lattice point $m\in\Delta\cap M$ and $a_m\in k$ is a general
choice of coefficient. So $z^0$ is the section which vanishes
once on each toric boundary component of $\PP_{\Delta}$.

Now for most choices of $\check h$, this does not define a toric degeneration:
the singularities are too bad. However, there is a natural toric variety
in which to describe this degeneration. 
Set
\[
\tilde\Delta:=\{(m,l)|m\in\Delta,\quad l\ge\check h'(m)\}\subseteq M_{\RR}
\oplus\RR.
\]
The assumption that $\check h'$ is convex implies $\tilde\Delta$ is
convex, but of course is non-compact.
Let $\tilde\Sigma$ be the normal fan to $\tilde\Delta$ in $N_{\RR}\oplus
\RR$. Let $X(\tilde\Sigma)$ denote
the toric variety defined by the fan $\tilde\Sigma$.
Then $\tilde\Delta$ is the Newton polytope of a line bundle $\shL$
on $X(\tilde\Sigma)$, and terms of the form $t^{\check h'(m)}z^m$
can be interpreted as sections of this line bundle, corresponding to
$(m,\check h'(m))\in \tilde\Delta\subseteq M_{\RR}\oplus\RR$. In addition, 
there is a natural map $X(\tilde\Sigma)\rightarrow\AA^1$ defined by projection
onto $\RR$, and this map defines the regular function $t$. 
Thus (\ref{toricdegenexam}) defines a hypersurface $\X_{\Delta}\subseteq
X(\tilde\Sigma)$ and $t$ defines a map $\X_{\Delta}\rightarrow\AA^1$.
This is a toric degeneration. 

Without going into much detail, choosing a star decomposition of
$\nabla$ as in \S 4 and a good polyhedral decomposition $\P$ of
$\partial\nabla^{\check h}$ is essentially the same as choosing
a refinement $\tilde\Sigma'$ of $\Sigma'$ with particularly nice
properties.
This yields a partial resolution $X(\tilde\Sigma')
\rightarrow X(\tilde\Sigma)$. Then the proper transform of $\X_{\Delta}$
in $X(\tilde\Sigma')$,
$\X'_{\Delta}$, yields another toric degeneration $\X'_{\Delta}\rightarrow
\AA^1$. This is necessary to get a toric degeneration whose general
fibre is a MPCP resolution of a hypersurface in a toric variety.
For proofs, see \cite{GBB}, where the construction is generalized to
complete intersections in toric varieties. \qed
\end{example}

Given a toric degeneration $f:\X\rightarrow D$, 
we can build the \emph{dual intersection
complex} $(B,\P)$ of $f$, as follows. Here $B$ is an integral affine
manifold with singularities, and $\P$ is a \emph{polyhedral decomposition}
of $B$, a decomposition of $B$ into lattice polytopes. In fact, we 
will construct $B$ as a union of lattice polytopes.
Specifically,
let the normalisation of $\X_0$, $\tilde \X_0$, be written as a disjoint
union $\coprod X_i$ of toric varieties $X_i$, $\nu:\tilde\X_0\rightarrow\X_0$
the normalisation. The {\it strata} of $\X_0$ are the elements of
the set
$$Strata(\X_0)=\{\nu(S)\,|\,\hbox{$S$ is a toric stratum of $X_i$ for some $i$}\}.$$
Here by toric stratum we mean the closure of a $(\CC^*)^n$ orbit.

Let $\{x\}\in Strata(\X_0)$ be a zero-dimensional stratum. 
Applying Definition \ref{toricdegen} (4) to a neighbourhood of $x$,
there is a toric variety $Y_x$ 
such that in
a neighbourhood of $x$, $f:\X\rightarrow D$ is locally isomorphic to
$f_x:Y_x\rightarrow\CC$, where $f_x$ is given by a monomial.
Now the condition that $f_x$ vanishes precisely
once along each toric divisor of $Y_x$ is the statement
that $Y_x$ is Gorenstein, and as such, it arises 
as in Exercise \ref{gorensteinlimit}. Indeed, let $M,N$ be as usual,
with $\rank M=\dim\X_0$. Then there
is a lattice polytope $\sigma_x\subseteq M_{\RR}$ such that $C(\sigma_x)
=\{(rm,r)|m\in\sigma, r\ge0\}$
is the cone defining the toric variety $Y_x$. As we saw in Exercise
\ref{gorensteinlimit}, a small neighbourhood of $x$ in $\X$ should contribute
a copy of $\sigma_x$ to $B$, which provides the motivation for our construction.
We can now describe how to construct $B$ by gluing together the polytopes
\[\{\sigma_x\,|\, \{x\}\in Strata(\X_0)\}.\]
We will do this in the case that every irreducible component
of $\X_0$ is in fact itself normal so that $\nu:X_i\rightarrow \nu(X_i)$ is an
isomorphism. The reader may be able to imagine the more
general construction.
With this normality assumption, there is a one-to-one inclusion reversing
correspondence between faces of
$\sigma_x$ and elements of $Strata(\X_0)$ containing $x$. We can then
identify faces of $\sigma_x$ and $\sigma_{x'}$ if they correspond
to the same strata of $\X_0$. Some argument is necessary to show that this
identification can be done via an integral affine transformation, but
again this is not difficult.

Making these identifications, one obtains $B$. One can then prove

\begin{lemma} If $\X_0$ is complex $n$ dimensional, then $B$ is an
real $n$ dimensional manifold.
\end{lemma}

See \cite{PartI}, Proposition 4.10 for a proof.

Now so far $B$ is just a topological manifold, constructed by gluing together
lattice polytopes. Let 
\[
\P=\{\sigma\subseteq B| \hbox{$\sigma$ is a face of $\sigma_x$ 
for some zero-dimensional stratum $x$}\}.
\]
There is a one-to-one inclusion reversing correspondence between strata
of $\X_0$ and elements of $\P$. 

It only remains to give $B$ an affine structure with singularities.
Let $\Bar(\P)$ be the first barycentric subdivision of $\P$, and let
$\Gamma$ be the union of simplices in $\Bar(\P)$ not containing
a vertex of $\P$ or intersecting the interior of a maximal cell of
$\P$. If we then set $B_0:=B\setminus\Gamma$, we can define an affine
structure on $B_0$ as follows.
$B_0$ has an open cover
\[
\{W_{\sigma}|\hbox{$\sigma\in\P$ maximal}\}\cup
\{W_v|\hbox{$v\in\P$ a vertex}\}
\]
where $W_{\sigma}=\Int(\sigma)$, the interior of $\sigma$, and 
\[
W_v=\bigcup_{\tau\in\Bar(\P)\atop v\in\tau}\Int(\tau)
\]
is the (open) star of $v$ in $\Bar(\P)$, just as in \S 4.

We now define charts.
As a maximal cell of $\P$
is of the form $\sigma_x\subseteq M_{\RR}$, this inclusion 
induces a natural affine chart $\psi_{\sigma}:\Int(\sigma)\rightarrow
M_{\RR}$.
On the other hand, a vertex $v$ of $\P$ corresponds to a codimension $0$
stratum of $\X_0$, i.e.\ to an irreducible component $X_i$ for some $i$.
Because this is a compact toric variety, it is defined by a complete
fan $\Sigma_i$ in $M_{\RR}$. Furthermore, there is a one-to-one
correspondence between $p$-dimensional
cones of $\Sigma_i$ and $p$-dimensional cells
of $\P$ containing $v$ as a vertex, as they both correspond to 
strata of $\X_0$ contained in $X_i$. There is then a continuous map
\[
\psi_v:W_v\rightarrow M_{\RR}
\]
which takes $W_v\cap\sigma$, for any $\sigma\in\P$ containing $v$ as a vertex,
into the corresponding cone of $\Sigma_i$
\emph{integral affine linearly}. Such a map is uniquely determined
by the combinatorial correspondence and the requirement that it
be integral affine linear on each cell. It is then obvious these
charts define an integral affine structure on $B_0$. Thus we have
constructed $(B,\P)$.

\begin{example} Let $f:\X\rightarrow D$ be a degeneration of elliptic
curves to an $I_n$ fibre. Then $B$ is the circle $\RR/n\ZZ$, decomposed
by $\P$ into $n$ line segments of length one.
\end{example}

\begin{example}
Continuing with Example \ref{toricdegenexample2}, the dual intersection
complex constructed from the toric degenerations $\X'_{\Delta}\rightarrow
\AA^1$ is the affine manifold with singularities structure
with polyhedral decomposition $(B,\P)$
constructed on $B=\partial\nabla^{\check h}$ in \S 4. This is not particularly
difficult to show. Again, for the proof and
more general complete intersection case,
see \cite{GBB}.
\end{example}

Is the dual intersection complex the right affine manifold with singularities?
The following theorem provides evidence for this, and gives the connection
between this construction and the SYZ conjecture.

\begin{theorem}
\label{complextheorem}
Let $\X\rightarrow D$ be a toric degeneration, with dual intersection
complex $(B,\P)$. Then there is an open set $U\subseteq B$ such that 
$B\setminus U$ retracts onto the discriminant locus $\Gamma$ of $B$,
such that $\X_t$ contains an open subset $\U_t$ which is isomorphic
as complex manifolds to a small deformation of a twist of $X_{\epsilon}(U)$,
where $\epsilon=O(-1/\ln|t|)$. 
\end{theorem}

We will not be precise here about what we mean by small deformation;
by twist, we mean a twist of the complex structure of $X_{\epsilon}(U)$
by a $B$-field. See \cite{Announce} 
for a much more precise statement; the above statement
is meant to give a feel for what is true. The proof, along with much more
precise statements, will eventually appear in \cite{tori}.

\medskip

If $\X\rightarrow D$ is a \emph{polarized} toric degeneration, i.e.\ if there
is a relatively ample line bundle $\shL$ on $\X$, then we can construct
another affine manifold with singularities and polyhedral decomposition
$(\check B,\check\P)$, which we call the \emph{intersection complex},
as follows.

For each irreducible component $X_i$ of $\X_0$, $\shL|_{X_i}$ is an ample
line bundle on a toric variety. Let $\sigma_i\subseteq N_{\RR}$
denote the Newton polytope of this line bundle. There is then
a one-to-one inclusion preserving correspondence between strata of
$\X_0$ contained in $X_i$ and faces of $\sigma_i$. We can then glue together
the $\sigma_i$'s in the obvious way: if $Y$ is a codimension one
stratum of $\X_0$, it is contained in two irreducible components $X_i$ and
$X_j$, and defines faces of $\sigma_i$ and $\sigma_j$. These faces are affine
isomorphic because they are both the Newton polytope of $\shL|_Y$, and
we can then identify them in the canonical way. Thus we obtain a topological
space $\check B$ with a polyhedral decomposition $\check\P$. We give
it an affine structure with singularities in a similar manner as before.
Again, let $\Gamma$ be the union of simplices in $\Bar(\check\P)$
not containing a vertex of $\check\P$ or intersecting the interior of
a maximal cell of $\check\P$. Setting $\check B_0:=\check B\setminus\Gamma$,
this again has an open cover
\[
\{W_{\sigma}|\hbox{$\sigma\in\check\P$ maximal}\}\cup\{W_v|\hbox{$v\in\check\P$
a vertex}\}.
\]
As usual, as $W_{\sigma}$ is the interior of $\sigma$, it comes along
with a canonical affine structure. On the other hand, a vertex $v$
of $\check\P$ corresponds to a dimension zero stratum $x$ of $\X_0$,
and associated to $x$ is the polytope $\sigma_x\subseteq M_{\RR}$.
Let $\check\Sigma_x$ be the normal fan to $\sigma_x$ in $N_{\RR}$. 
Then there is
a one-to-one inclusion preserving correspondence between cones in
$\check\Sigma_x$ and strata of $\X_0$ containing $x$. This correspondence
allows us to define a chart
\[
\check\psi_v:W_v\rightarrow N_{\RR}
\]
which takes $W_v\cap\check\sigma$, 
for any $\check\sigma\in\check\P$ containing $v$ as a vertex,
into the corresponding cone of $\check\Sigma_x$
in an integral affine linear way. 
This gives the manifestly integral affine structure on $\check B_0$,
and hence defines the intersection complex $(\check B,\check\P)$.

Analogously to Theorem \ref{complextheorem}, we expect

\begin{conjecture}
Let $\X\rightarrow D$ be a polarized toric degeneration, 
with intersection complex
$(\check B,\check\P)$. 
Let $\omega_t$ be a K\"ahler form on $\X_t$ representing the first Chern
class of the polarization.
Then there is an open set $\check U\subseteq \check B$ 
such that $\check B\setminus \check U$ retracts onto the discriminant locus
$\Gamma$ of $\check B$, such that $\X_t$ is a symplectic compactification
of $\check X(\check U)$ for any $t$.
\end{conjecture}

I don't expect this to be particularly difficult: it should be amenable
to the techniques of Ruan \cite{RuanJSG}, but has not been carried out.

\medskip

The relationship between the intersection complex and the dual intersection
complex can be made more precise by introducing multi-valued piecewise
linear functions, in analogy with the multi-valued convex functions of
Definition \ref{multivaluedconvex}:

\begin{definition}
Let $B$ be an affine manifold with singularities with polyhedral
decomposition $\P$. Then a multi-valued piecewise linear function
$\varphi$ on $B$ is a collection of continuous functions on an open cover
$\{(U_i,\varphi_i)\}$ such that $\varphi_i$ is affine linear on each
cell of $\P$ intersecting $U_i$, and on $U_i\cap U_j$, $\varphi_i-\varphi_j$
is affine linear. Furthermore, for any $\sigma\in\P$,
in a neighbourhood of each point $x\in U_i\cap
Int(\sigma)$, there is an affine linear function $\psi$ such that 
$\varphi_i-\psi$ is zero on $\sigma$.
\end{definition}

To explain this last condition, and to clarify additional
structure on $B$, let us examine a property of the polyhedral decomposition
$\P$ of $B$ when $(B,\P)$ is a dual intersection complex.

Consider any $p$-dimensional cell $\sigma\in\P$. This corresponds to
an $n-p$-dimensional stratum $X_{\sigma}\subseteq\X_0$, and as such,
it is a toric variety defined by a fan $\Sigma_{\sigma}$ in $\RR^{n-p}$.
Now for any vertex $v$ of $\sigma$, $X_{\sigma}$ is a toric stratum
of the irreducible component $X_v$ of $\X_0$. Thus $\Sigma_{\sigma}$ can
be obtained as a \emph{quotient fan} of $\Sigma_v$. In other words, there
is a $p$-dimensional cone $K_{\sigma}$ of $\Sigma_v$ corresponding
to $\sigma$ such that
\[
\Sigma_{\sigma}=\Sigma_v(K_{\sigma}):=\{(K+\RR K_{\sigma})/\RR K_{\sigma}|
K\in\Sigma_v,K\supseteq K_{\sigma}\}.
\]
In particular, there is an open neighbourhood $U_{v,\sigma}$
of $\Int(K_{\sigma})$ and an integral linear map $S_{v,\sigma}:U_{v,\sigma}
\rightarrow \RR^{n-p}$ such that
\[
\{S_{v,\sigma}^{-1}(K)| K\in \Sigma_{\sigma}\}=
\{U_{v,\sigma}\cap K| K\in\Sigma_v, K\supseteq K_{\sigma}\}.
\]
Let $U_{\sigma}$ be a small open neighbourhood of $\Int(\sigma)$
in $B$;
if taken sufficiently small the maps $S_{v,\sigma}$ can be viewed
as being defined on open subsets of $U_{\sigma}$ and patch to give an integral
affine submersion $S_{\sigma}:U_{\sigma}\rightarrow\RR^{n-p}$, where
$U_{\sigma}=\bigcup_{v\in\sigma} U_{v,\sigma}$ is an open neighbourhood
of $\Int(\sigma)$. This map has the 
property that 
\[
\{S_{\sigma}^{-1}(K)|K\in\Sigma_{\sigma}\}=\{U_{\sigma}\cap\tau|
\tau\supseteq\sigma,\tau\in\P\}.
\]
In general we call a polyhedral decomposition
\emph{toric} if for all $\sigma\in\P$ there is always such an integral affine
linear map 
$S_{\sigma}:U_{\sigma}\rightarrow
\RR^{n-p}$ from a neighbourhood $U_{\sigma}$ of $\Int(\sigma)$
and a fan $\Sigma_{\sigma}$ in $\RR^p$ with the above property.
(See \cite{PartI}, Definition 1.22 for a perhaps too precise definition
of toric polyhedral decompositions. The definition there is complicated
by allowing cells to be self-intersecting, or equivalently, allowing
irreducible components of $\X_0$ to be non-normal.)
We can think of the fan $\Sigma_{\sigma}$ as being the fan structure of
$\P$ transverse to $\sigma$ at a point in the interior of $\sigma$.
The main point for a dual intersection complex
is that this fan structure is determined by
$X_{\sigma}$, and this is independent of the choice of the point
in the interior of $\sigma$. 

Now let us return to piecewise linear functions. Suppose we are given a 
polarized toric degeneration $\X\rightarrow D$. We in fact obtain a piecewise
linear function $\varphi$ on the dual intersection complex $(B,\P)$
as follows. Restricting to any toric stratum $X_{\sigma}$, $\shL|_{X_{\sigma}}$
is determined completely by an integral piecewise linear function
$\bar\varphi_{\sigma}$ on $\Sigma_{\sigma}$, well-defined up to
a choice of linear function. Pulling back this piecewise linear function
via $S_{\sigma}$ to $U_{\sigma}$, we obtain a collection of piecewise
linear functions $\{(U_{\sigma},\varphi_{\sigma})|\sigma\in\P\}$.
The fact that $(\shL|_{X_{\tau}})|_{X_{\sigma}}=\shL|_{X_{\sigma}}$
for $\tau\subseteq\sigma$ implies that on overlaps $\varphi_{\sigma}$
and $\varphi_{\tau}$ differ by at most a linear function. So $\{(U_{\sigma},
\varphi_{\sigma})\}$ defines a multi-valued piecewise linear function.
The last condition in the definition of multi-valued piecewise linear
function then reflects the need for the function to be locally
a pull-back of a function via $S_{\sigma}$ in a neighbourhood of $\sigma$.
In fact, given any multi-valued piecewise linear function $\varphi$
on $(B,\P)$ with $\P$ a toric polyhedral decomposition of $B$, $\varphi$ 
is determined by functions $\bar\varphi_{\sigma}$ on $\Sigma_{\sigma}$
for $\sigma\in\P$, via pull-back by $S_{\sigma}$.

If $\shL$ is ample, then the piecewise linear function determined by
$\shL|_{X_{\sigma}}$ is strictly convex. So we say a multi-valued piecewise
linear function is \emph{strictly convex} if $\bar\varphi_{\sigma}$ is
strictly convex for each $\sigma\in\P$. 

\medskip

Now suppose we are given abstractly a triple $(B,\P,\varphi)$ with $B$
an integral affine manifold with singularities with a toric polyhedral
decomposition $\P$, and $\varphi$ a strictly convex multi-valued
piecewise linear function on $B$. Then we construct the \emph{discrete
Legendre transform} $(\check B,\check\P,\check\varphi)$ of $(B,\P,\varphi)$
as follows.

$\check B$ will be constructed by gluing together Newton polytopes.
If we view, for $v$ a vertex of $\P$, the fan $\Sigma_v$ as living
in $M_{\RR}$, then the Newton polytope of $\bar\varphi_v$ is
\[
\check v=\{x\in N_{\RR}|\langle x,y\rangle\ge-\bar\varphi_v(y)
\quad\forall y\in M_{\RR}\}.
\]
There is a one-to-one order reversing correspondence between faces of
$\check v$ and cells of $\P$ containing $v$. Furthermore, if
$\sigma$ is the smallest cell of $\P$ containing two vertices $v$ and $v'$,
then the corresponding faces of $\check v$ and $\check v'$ are integral 
affine isomorphic, as they are both isomorphic to the Newton polytope
of $\bar\varphi_{\sigma}$. Thus we can glue $\check v$ and $\check v'$
along this common face. After making all these identifications, we obtain a cell
complex $(\check B,\check\P)$, which is really just the dual cell complex
of $(B,\P)$. Of course, we have some additional information, namely an
affine structure on the interior of each maximal cell of $\check\P$.
To give $\check B$ an integral affine structure with singularities,
one proceeds as usual, using this affine structure along with
an identification of a neighbourhood of each vertex of $\check\P$ with
the normal fan of the corresponding maximal cell of $\P$.

Finally, the function $\varphi$ has a discrete Legendre transform
$\check\varphi$ on $(\check B,\check\P)$. We have no choice but to
define $\check\varphi$ in a neighbourhood of a vertex $\check\sigma\in
\check\P$ dual to a maximal cell $\sigma\in\P$ to be a piecewise
linear function whose Newton polytope is $\sigma$, i.e.\
\[
\overline{\check\varphi}_{\check\sigma}(y)
=-\inf\{\langle y,x\rangle| x\in\sigma\subseteq M_{\RR}\}.
\]
This gives $(\check B,\check\P,\check\varphi)$, the discrete
Legendre transform of $(B,\P,\varphi)$. If $B$ is $\RR^n$, then
this coincides with the classical notion of a discrete Legendre
transform. The discrete Legendre transform has several relevant 
properties:
\begin{itemize}
\item The discrete Legendre transform of $(\check B,\check\P,\check\varphi)$
is $(B,\P,\varphi)$.
\item If we view the underlying topological spaces $B$ and $\check B$
as being identified by being the underlying space of dual cell complexes,
then $\Lambda_{B_0}\cong \check\Lambda_{\check B_0}$ and
$\check\Lambda_{B_0}\cong\Lambda_{\check B_0}$, where the subscript
denotes which affine structure is being used to define $\Lambda$ or
$\check\Lambda$. 
\end{itemize}

This hopefully makes it clear that the discrete Legendre transform
is a suitable replacement for the duality provided to us by the
Legendre transform of \S 2.

Finally, it leads to what we may think of as an \emph{algebro-geometric
SYZ procedure}. In analogy with the procedure suggested in \S 5, we
follow these steps:
\begin{enumerate}
\item We begin with a toric degeneration of Calabi-Yau
manifolds $\X\rightarrow D$ with an ample polarization.
\item Construct $(B,\P,\varphi)$ from this data, as explained above.
\item Perform the discrete Legendre transform to obtain $(\check B,
\check\P,\check\varphi)$.
\item Try to construct a polarized degeneration of Calabi-Yau
manifolds $\check\X\rightarrow D$ whose dual intersection
complex is $(\check B,\check\P,\check\varphi)$.
\end{enumerate}

\begin{example} The discrete Legendre transform enables us to reproduce 
Baty\-rev duality. Returning to the construction of \S 4 and
Example \ref{toricdegenexample2}, choosing
a strictly convex piecewise linear function on $\tilde\Sigma'$ corresponding
to a line bundle $\shL$ induces
a polarization of $\X'_{\Delta}$. This then gives us a strictly
convex multi-valued piecewise linear function $\varphi_{\shL}$
on $(B,\P)$, hence
a discrete Legendre transform $(\check B,\check\P,\check\varphi_{\shL})$.
In \cite{GBB} I showed that this is the dual intersection complex associated
to some choice of subdivision $\widetilde{\check\Sigma'}$ of $\widetilde
{\check\Sigma}$ obtained by interchanging the roles of $\nabla$ and
$\Delta$ in the construction of \S 4. As an exercise, you can check 
the following for
yourself. If we take $\check h=\check\varphi$,
and in addition define 
\[
h=\varphi:N_{\RR}\rightarrow \RR
\]
to take the value $1$ on the primitive generator of each one-dimensional
cone on $\Sigma$, the normal fan to $\Delta$, then from \S 4 we obtain
an affine structure with singularities on
$B=\partial\nabla$, and completely symetrically using $h$
we also obtain such a structure on $\check B=\partial\Delta$.
These manifolds come with
polyhedral decomposition $\P$ and $\check\P$ consisting of all proper faces
of $\nabla$ and $\Delta$ respectively.
The anti-canonical polarizations on $\PP_{\Delta}$ and $\PP_{\nabla}$
induce multi-valued piecewise linear
functions $\psi,\check\psi$ on  $B$ and $\check B$ respectively.
Then show $(B,\P,\psi)$ and $(\check B,\check\P,\check\psi)$ are
discrete Legendre transforms of each other.

Thus Batyrev (and Batyrev-Borisov) duality is a special case of this
construction.
\end{example}

The only step missing in this mirror symmetry algorithm is the last: 
\begin{question}[The reconstruction problem, Version II]
\label{reconstruct2}
Given $(B,\P,\varphi)$, is it possible to construct a polarized toric
degeneration $\X\rightarrow D$ whose dual intersection complex is
$(B,\P,\varphi)$?
\end{question}

It is fairly obvious how to reconstruct the central fibre $\X_0$ of
a degeneration from the data $(\check B,\check\P,\check\varphi)$, and
we will see this explicitly in \S 8.
One could naively hope that this reducible variety has good deformation
theory and it can be smoothed. However, in general its deformation theory
is ill-behaved.

As initially observed in the normal crossings case by Kawamata and
Namikawa in \cite{KN}, one needs to put some additional structure on $\X_0$
before it has good deformation theory. This structure is a \emph{log
structure}, and introducing log structures allows us to study
many aspects of mirror symmetry directly on the degenerate fibre itself.

We shall do this in the next section, but first, let me address the
question of how general this mirror symmetry construction might be:

\begin{conjecture} If $f:\X\rightarrow D$ is a large complex
structure limit degeneration, then $f$ is birationally equivalent
to a toric degeneration $f':\X'\rightarrow D$.
\end{conjecture}

The condition of being a large complex structure limit, as defined
by Morrison in \cite{Morr}, is a stronger one than maximally unipotent.
Why should I imagine
something like this to be true? Well, fantasizing freely, we would expect
that after choosing a polarization and Ricci-flat metric on fibres $\X_t$,
we have a sequence $(\X_t,g_t)$ converging to $B$ an affine manifold with
singularities. Now in general an affine manifold with singularities need
not arise as the dual intersection complex of a toric degeneration,
first of all because it need not have a toric polyhedral decomposition. For
example, even in two dimensions there are orbifold singularities (corresponding
to singular elliptic fibres which are not semi-stable) which do not
arise in dual intersection complexes of toric degenerations, yet can
occur as the base of a special Lagrangian fibration on a K3 surface.
However, the \emph{general} base does arise as the dual intersection complex
of a toric degeneration in the K3 case. The hope is that the condition
of large complex structure limit forces the singularities of $B$ to
be ``sufficiently general'' so that one can construct a nice toric polyhedral
decomposition $\P$ on $B$, and from this construct a toric degeneration.
Presumably, this toric degeneration will be, if picked correctly, 
birational to the original one. 

This argument of course is rather hand-wavy, but I believe it provides
some moral expectation that there might be a large class of degenerations
for which our method applies. I note that one 
can prove the conjecture in the case of K3 surfaces.

We now come to the technical heart of the program laid out in \cite{PartI}.
Some aspects of this program are quite technical, so the goal here is
to explain the highlights of \cite{PartI} as simply possible.

\section{Log structures}

We first introduce the log structures of Fontaine-Illusie and Kato
(\cite{Illu}, \cite{K.Kato}). 

\begin{definition}
A log structure on a scheme (or analytic space) $X$ is a (unital) homomorphism
$$\alpha_X:\shM_X\rightarrow \O_X$$
of sheaves of (multiplicative and commutative) monoids inducing an isomorphism
$\alpha_X^{-1}(\O_X^{\times})\rightarrow \O_X^{\times}$. The
triple $(X,\shM_X,\alpha_X)$ is then called a {\it log space}.
We often write the whole package as $X^{\dagger}$.
\end{definition}

A morphism of log spaces $F:X^{\dagger}\rightarrow Y^{\dagger}$ consists
of a morphism $\underline{F}:X\rightarrow Y$ of underlying
spaces together with a homomorphism $F^{\#}:\underline{F}^{-1}(\shM_Y)
\rightarrow\shM_X$ commuting with the structure homomorphisms:
$$\alpha_X\circ F^{\#}=\underline{F}^*\circ\alpha_Y.$$

The key examples:

\begin{examples}
\label{logexamples}
(1) Let $X$ be a scheme and $Y\subseteq X$ a closed subset of
codimension one. Denote by $j:X\setminus Y\rightarrow X$
the inclusion. Then the inclusion
$$\alpha_X:\shM_X=j_*(\O_{X\setminus Y}^{\times})\cap\O_X\rightarrow
\O_X$$
of the sheaf of regular functions with zeroes contained in $Y$ is a log
structure on $X$. This is called a \emph{divisorial log structure} on
$X$.

(2) A {\it prelog structure}, i.e.\ an arbitrary homomorphism of
sheaves of monoids $\varphi:\shP\rightarrow\O_X$, defines
an associated log structure $\shM_X$ by
$$\shM_X=(\shP\oplus\O_X^{\times})/\{(p,\varphi(p)^{-1})|p\in
\varphi^{-1}(\O_X^{\times})\}$$
and $\alpha_X(p,h)=h\cdot\varphi(p)$.

(3) If $f:X\rightarrow Y$ is a morphism of schemes and $\alpha_Y:\shM_Y
\rightarrow\O_Y$ is a log structure on $Y$, then the prelog structure
$f^{-1}(\shM_Y)\rightarrow\O_X$ defines an associated log structure
on $X$, the {\it pull-back log structure}.

(4) In (1) we can pull back the log structure on $X$ to $Y$ using
(3). Thus in particular, if $\X\rightarrow D$ is a toric
degeneration, the inclusion $\X_0\subseteq\X$ gives a log
structure on $\X$ and an induced log structure on $\X_0$. Similarly
the inclusion $0\in D$ gives a log structure on $D$ and
an induced one on $0$. Here $\M_0=\CC^{\times}\oplus\NN$,
where $\NN$ is the (additive) monoid of natural (non-negative) numbers,
and
$$\alpha_0(h,n)=\begin{cases}h& n=0\\ 0&n\not=0.\end{cases}$$
$0^{\dagger}$ is usually called the standard log point.

We then have log morphisms $\X^{\dagger}\rightarrow D^{\dagger}$ and
$\X_0^{\dagger}\rightarrow 0^{\dagger}$.

(5) If $\sigma\subseteq M_{\RR}=\RR^n$ is a
convex rational polyhedral cone, $\dual{\sigma}\subseteq
N_{\RR}$ the dual cone, let $P=\dual{\sigma}\cap N$: this is a monoid.
The affine toric variety defined by $\sigma$ can be written as
$X=\Spec \CC[P]$. 
We then have a pre-log structure induced by the homomorphism of
monoids
$$P\rightarrow \CC[P]$$
given by $p\mapsto z^p$. There is then an associated log
structure on $X$. This is in fact the same as the log structure
induced by $\partial X\subseteq X$, where $\partial X$
is the toric boundary of $X$, i.e.\ the union of toric divisors of $X$.

If $p\in P$, then the monomial $z^p$ defines a map
$f:X\rightarrow \Spec \CC[\NN]\quad (=\Spec \CC[t])$ which is a log morphism
with the log structure on $\Spec \CC[\NN]$ induced similarly by
$\NN\rightarrow\CC[\NN]$.
The fibre $X_0=\Spec \CC[P]/(z^p)$ is a subscheme of $X$,
and there is an induced log structure on $X_0$, and a map $X_0^{\dagger}
\rightarrow 0^{\dagger}$ as in (4). $f$ is an example of a 
\emph{log smooth} morphism. Essentially all log smooth morphisms
are \'etale locally of this form (if $\NN$ is replaced by a more general
monoid). See \cite{F.Kato} for details.

Condition (4) of Definition \ref{toricdegen} in fact implies
that locally, away from $Z$, $\X^{\dagger}$ and $\X_0^{\dagger}$ are
of the above form. So we should view $\X^{\dagger}\rightarrow D^{\dagger}$
as log smooth away from $Z$, and from the log point of view, $\X_0^{\dagger}$
can be treated much like a non-singular scheme away from $X$. We will
see this explicitly below when we talk about differentials. \qed
\end{examples}

\medskip

On a log scheme $X^{\dagger}$ there is always an exact sequence
\[
1\mapright{} \O_X^{\times}\mapright{\alpha^{-1}}\M_X\mapright{}
\overline{\M}_X\mapright{}0,
\]
where we write the quotient sheaf of monoids $\overline{\M}_X$
additively. We call $\overline{\M}_X$ the \emph{ghost sheaf}
of the log structure. I like to view $\overline{\M}_X$ as specifying
the combinatorial information associated to the log structure. For 
example, if $X^{\dagger}$ is induced by the Cartier divisor $Y\subseteq
X$ with $X$ normal, 
then the stalk $\overline{\M}_{X,x}$ at $x\in X$ is the monoid 
of effective Cartier divisors on a neighbourhood of $x$ supported
on $Y$.

\begin{xca}
\label{ghostexercise}
Show that in Example \ref{logexamples}, (5), $\overline{\M}_{X,x}=
P$ if $\dim\sigma=n$ and $x$ is the unique zero-dimensional torus orbit
of $X$. More generally, 
\[
\overline{\M}_{X,x}={\dual{\tau}\cap N\over \tau^{\perp}\cap N}
=\Hom_{monoid}(\tau\cap M,\NN),
\]
when $x\in X$ is in the torus orbit 
corresponding to a face $\tau$ of $\sigma$. In particular, $\tau$
can be recovered as $\Hom_{monoid}(\overline{\M}_{X,x},\RR_{\ge 0}^+)$,
where $\RR_{\ge 0}^+$ is the additive monoid of non-negative real
numbers.
\qed
\end{xca}

Another important fact for us is that if $f:Y\rightarrow X$ is a morphism
with $X$ carrying a log structure, and $Y$ is given the pull-back log
structure, then $\overline{\M}_Y=f^{-1}\overline{\M}_X$. In the case
that $\M_X$ is induced by an inclusion
of $Y\subseteq X$, $\overline{\M}_X$ is supported on $Y$, so we can
equate $\overline{\M}_X$ and $\overline{\M}_Y$, the ghost sheaves
for the divisorial log structure on $X$ and its restriction to $Y$.
Putting this together with Exercise \ref{ghostexercise}
and the definition of dual intersection complex, we see that given
a toric degeneration $\X\rightarrow D$ the dual intersection complex
\emph{completely determines} the ghost sheaf $\overline{\M}_{\X}=
\overline{\M}_{\X_0}$ off of $Z$. We in fact take the view that
anyway the log structure on $Z$ is not particularly well-behaved,
and we always ignore it on $Z$. In fact, given a log structure 
$\M_{\X_0\setminus
Z}$ on $\X_0\setminus Z$, this defines a push-forward log structure
$\M_{\X_0}:=j_*\M_{\X_0\setminus Z}$. There is an induced map
$\alpha:\M_{\X_0}\rightarrow\O_{\X_0}$, as $j_*\O_{\X_0\setminus Z}
=\O_{\X_0}$ because $\X_0$ is Cohen-Macaulay. Thus in what follows, if we
have determined a log structure on $\X_0\setminus Z$,
we just as well get a log structure on $\X_0$ and will not
concern ourselves with the behaviour of this log structure along $Z$.

All this gives the necessary hint for working backwards, to go from
$(B,\P)$ to $\X_0^{\dagger}$. Suppose we are given an integral
affine manifold with singularities $B$ with \emph{toric} 
polyhedral decomposition
$\P$. At each vertex $v$ of $\P$, $\P$ locally looks like a fan
$\Sigma_v$, defining a toric variety $X_v$. For every edge $\omega
\in\P$ with endpoints $v$ and $w$, $\omega$ defines a ray in both fans
$\Sigma_v$ and $\Sigma_w$, hence toric divisors $D^v_{\omega}
\subseteq X_v$, $D^w_{\omega}\subseteq X_w$. The condition that $\P$
is a toric polyhedral decomposition tells us that $D^v_{\omega}$
and $D^w_{\omega}$ are isomorphic toric varieties, and we can choose a 
torus equivariant isomorphism $s_{\omega}:D^v_{\omega}\rightarrow D^w_{\omega}$
for each edge $\omega$. If we choose these gluing maps to satisfy a
certain compatibility condition on codimension two strata (we leave it
to the reader to write down this simple compatibility condition), then
we can glue together the $X_v$'s to obtain, in general, an algebraic
space we write as $X_0(B,\P,s)$, where $s=(s_{\omega})$ is the collection
of gluing maps. (In \cite{PartI}, we describe the gluing data in a slightly
different, but equivalent, way). We call $s$ \emph{closed gluing data}.
This is how we construct a potential central fibre of a toric degeneration.

Now $X_0(B,\P,s)$ cannot be a central fibre of a toric degeneration unless
it carries a log structure of the correct sort. There are many reasons
this may not happen. First, if $s$ is poorly chosen, there may be
zero-dimensional strata of $X_0(B,\P,s)$ which do not have neighbourhoods
locally \'etale isomorphic to the toric boundary of an affine toric variety;
this is a minimum prerequisite. As a result, we have to
restrict attention to closed
gluing data induced by what we call \emph{open gluing data}. 
Explicitly, each maximal cell $\sigma\in\P$ defines an affine toric
variety $U(\sigma)$ given by the cone $C(\sigma)\subseteq M_{\RR}\oplus
\RR$, assuming we view $\sigma\subseteq M_{\RR}$ as a lattice polytope.
Let $V(\sigma)\subseteq U(\sigma)$ be the toric boundary. It turns out,
as we show in \cite{PartI}, that a necessary condition for $X_0(B,\P,s)$
to be the central fibre of a toric degeneration is that it is obtained
by dividing out $\coprod_{\sigma\in\P_{\max}} V(\sigma)$ by an
\'etale equivalence relation. 
In other words, we are gluing together the $V(\sigma)$'s to obtain an
algebraic space, and
those \'etale equivalence relations
which produce algebraic spaces of the form $X_0(B,\P,s)$ are easily
determined. This is carried out in detail in \cite{PartI}, \S 2. 
The construction there appears technically difficult because of the necessity
of dealing with algebraic spaces, but is basically straightforward.
The basic point is that if $\sigma_1,\sigma_2\in\P$ are two
maximal cells, with $\sigma_1\cap\sigma_2=\tau$, then $\tau$
determines faces of the cones $C(\sigma_1)$ and $C(\sigma_2)$,
hence open subsets $U_i(\tau)\subseteq U(\sigma_i)$, with
toric boundaries $V_i(\tau)\subseteq V(\sigma_i)$. Now in general
there is no \emph{natural} isomorphism between
$U_1(\tau)$ and $U_2(\tau)$: this
is a problem when $\sigma_1\cap\sigma_2\cap\Gamma\not=\emptyset$, where
$\Gamma$ is as usual the singular locus of $B$. However, crucially $V_1(\tau)$
and $V_2(\tau)$ are naturally 
isomorphic, and we can choose compatible equivariant
isomorphisms to obtain \emph{open gluing data}. Choosing open
gluing data allows us to define the \'etale equivalence relation:
we are just gluing any two sets $V(\sigma_1),V(\sigma_2)$
via the chosen isomorphism between $V_1(\tau)$ and $V_2(\tau)$.
Any choice of open gluing data $s$ gives rise in this way to an
algebraic space $X_0(B,\P,s)$, and to any choice of open gluing data
there is associated closed gluing data $s'$ such that
$X_0(B,\P,s)\cong X_0(B,\P,s')$.

The advantage of using open gluing data is that each $V(\sigma)$
for $\sigma\in\P_{\max}$ carries a log structure induced by the
divisorial log structure $V(\sigma)\subseteq U(\sigma)$. Unfortunately,
these log structures are not identified under the open gluing maps, precisely
because of a lack of a natural isomorphism between the $U_i(\tau)$'s
cited above. However,
the ghost sheaves of the log structures are isomorphic. 
So the ghost sheaves $\overline{\M}_{V(\sigma)}$
glue to give a ghost sheaf of monoids $\overline{\M}_{X_0(B,\P,s)}$.
Summarizing what we have said so far: (this is a combination
of results of \cite{PartI}, \S\S 2,4)

\begin{theorem}
Given $(B,\P)$, if $s$ is closed gluing data, and $\X_0=X_0(B,\P,s)$
is the central fibre of a toric degeneration $\X\rightarrow D$
with dual intersection complex $(B,\P)$, then $s$ is induced by open
gluing data and $\overline{\M}_{\X_0}|_{\X_0\setminus Z}
\cong\overline{\M}_{X_0(B,\P,s)}|_{\X_0\setminus Z}$.
\end{theorem}

This is as far as we can get with the combinatorics. The next point is
to attempt to construct $\M_{X_0(B,\P,s)}$. The idea is that
$\M_{X_0(B,\P,s)}$ is an extension of $\overline{\M}_{X_0(B,\P,s)}$
by $\O_{X_0(B,\P,s)}^{\times}$, so we are looking for some subsheaf of the sheaf
\[
\shExt^1(\overline{\M}_{X_0(B,\P,s)}^{\gp},\O_{X_0(B,\P,s)}^{\times}).
\]
Here the superscript $\gp$ denotes the Grothendieck group of the monoid.
Any extension of $\overline{\M}_{X_0(B,\P,s)}^{\gp}$ by 
$\O_{X_0(B,\P,s)}^{\times}$
gives rise to a sheaf of groups $\M_{X_0(B,\P,s)}^{\gp}$ surjecting
onto $\overline{\M}_{X_0(B,\P,s)}^{\gp}$, and the inverse image of
$\overline{\M}_{X_0(B,\P,s)}\subseteq\overline{\M}_{X_0(B,\P,s)}^{\gp}$
is a sheaf of monoids $\M_{X_0(B,\P,s)}$. Of course, one also needs
a map $\alpha:\M_{X_0(B,\P,s)}\rightarrow\O_{X_0(B,\P,s)}$,
and this complicates things a bit. To make a long story short, we
can identify a subsheaf of extensions which yield genuine log structures.
A section of this subsheaf determines a log structure on $X_0(B,\P,s)$
with the correct ghost sheaf. However, this is not precisely what we
want. What we really want is a log structure
on $X_0(B,\P,s)$ along with a log morphism $X_0(B,\P,s)^{\dagger}
\rightarrow 0^{\dagger}$ which is log smooth. (We will address
the question of the bad set $Z\subseteq\X_0$ shortly.)
We call such a structure a \emph{log smooth structure} on 
$X_0(B,\P,s)$. It turns out these structures
are given by certain sections of 
$\shExt^1(\overline{\M}^{\gp}_{X_0(B,\P,s)}/
\bar\rho,\O_X^{\times})$, 
where $\bar\rho$ is the canonical section of
$\overline{\M}_{X_0(B,\P,s)}$ whose germ at $\Mbar_{X_0(B,\P,s),\eta}=
\NN$ is $1$ for $\eta$ a generic point of an irreducible component
of $X_0(B,\P,s)$.
So in fact, we can identify a subsheaf of
$\shExt^1(\overline{\M}^{\gp}_{X_0(B,\P,s)}/\bar\rho,\O_X^{\times})$, 
which we call $\shLS_{X_0(B,\P,s)}$, whose sections determine a log structure
on $X_0(B,\P,s)$ \emph{and} a log smooth morphism 
$X_0(B,\P,s)^{\dagger}\rightarrow 0^{\dagger}$, i.e.\ a log
smooth structure.

The technical heart of \cite{PartI} is an explicit calculation of the
sheaf $\shLS_{X_0(B,\P,s)}$. This is carried out locally in
\cite{PartI}, Theorem 3.22, where the sheaf is calculated on the (\'etale)
open subsets $V(\sigma)$ of $X_0(B,\P,s)$, and globally in \cite{PartI}, Theorem
3.24. I will not state the precise results, but go into detail in 
a special case, which illustrates the most important features of the theory.

\begin{example} Suppose $X_0(B,\P,s)$ is normal crossings, i.e.\ every
cell of $\P$ is affine isomorphic to a standard simplex. Then we have
the local $\T^1$ sheaf, 
\[
\T^1=\shExt^1_{X_0(B,\P,s)}(
\Omega^1_{X_0(B,\P,s)/k},\O_{X_0(B,\P,s)}).
\]
This is a line bundle on
$S=Sing(X_0(B,\P,s))$. Then one can show $\shLS_{X_0(B,\P,s)}$ is the
$\O_S^{\times}$-torsor associated to $\T^1$.

This brings us back to Friedman's condition of $d$-semistability
\cite{Friedman}. A variety with normal crossings is \emph{$d$-semistable}
if $\T^1\cong\O_S$. Thus we recover Kawamata and Namikawa's result
\cite{KN}
showing that $X_0(B,\P,s)$ carries a normal crossings log structure
over $0^{\dagger}$ if and only if $X_0(B,\P,s)$ is $d$-semistable.
This is because, of course, the $\O_S^{\times}$-torsor associated
to $\T^1$ has a section if and only if $\T^1\cong\O_S$.

Now Theorem 3.24 of \cite{PartI} tells us that in general $\shLS_{X_0(B,\P,s)}$
is not a trivial $\O_S^{\times}$-torsor. The sheaf depends continuously
on $s$, but discretely on monodromy of the singularities of $B$.

Let's explain the latter point explicitly if $\dim B=2$. The irreducible
components of $S$ are in one-to-one correspondence with one-dimensional
cells of $\P$. If $\omega\in\P$ is such an edge, suppose it contains
one singularity of $B$ such that $\Lambda$ has monodromy
$\begin{pmatrix}1&n\\ 0&1\end{pmatrix}$ 
in a suitable basis around a loop around 
the singularity. Then $\T^1$ restricted to the one-dimensional stratum
$X_{\omega}\cong\PP^1$ of $X_0(B,\P,s)$ is $\O_{\PP^1}(n)$.

To make this statement completely accurate, one needs to define $n$
so that it is independent of the choice of basis and loop. To do this,
one chooses a loop which is counterclockwise with respect
to the orientation determined by the chosen basis of $\Lambda_b=\T_{B,b}$, 
where $b\in B$ is the base-point of the loop. 

If all the $n$'s appearing are positive, then for some choices of
gluing data $s$, we may hope to have a section $t$ of $\T^1$ which
vanishes only at a finite set of points $Z$. If $Z$ does not contain
a toric stratum (i.e.\ a triple point) then we obtain a log structure
on $X_0(B,\P,s)\setminus Z$ of the desired sort, hence a log
structure on $X_0(B,\P,s)$ (log smooth off of $Z$)
by push-forward. We then have

\begin{proposition}
In the situation of this example, with $\dim B=2$ and 
\[
t\in\Gamma(X_0(B,\P,s),\T^1)
\]
a section vanishing on a finite set $Z$ not containing a triple point,
there exists a smoothing $\X\rightarrow D$ of $X_0(B,\P,s)$
such that the singular locus of $\X$ is $Z\subseteq
\X_0=X_0(B,\P,s)$, and the induced log morphism $\X_0^{\dagger}
\rightarrow 0^{\dagger}$
coincides with $X_0(B,\P,s)^{\dagger}\rightarrow 0^{\dagger}$ 
determined by $t$.
\end{proposition}

The proof of this is a rather simple application of Friedman's or
Kawamata and Namikawa's results. To apply these results, however,
we need to deal with the singular set $Z$. This is done by normalizing
$X_0(B,\P,s)$, choosing to blow up one point in the inverse image
of each point of $Z$, and then regluing along the proper transform
of the conductor locus. This produces a $d$-semistable variety,
in the language of Friedman, or a log smooth scheme, which can
then be smoothed. (Such an approach seems difficult in higher dimensions.) 

On the other hand, if $n<0$ for some singular point of $B$, we run into
problems, and there is in fact no smoothing of $X_0(B,\P,s)$. This should
not be surprising for the following reason. If $n=-1$, it turns out
we would have to compactify the torus fibration
$X(B_0)$ by adding a strange sort of $I_1$
fibre over such a singular point. An $I_1$ fibre is an immersed
sphere, and the intersection multiplicity of the two sheets at 
the singular point of the fibre is $+1$ for an ordinary
$I_1$ fibre. However, when the monodromy is given by $n=-1$, the
intersection multiplicity is $-1$. This does not occur for a special
Lagrangian $T^2$-fibration, so it is not surprising we can't construct
a smoothing in this case. \qed
\end{example}

If $\dim B=2$ and $n>0$ for all singularities on $B$, then we say
$B$ is \emph{positive}. One can generalize this notion of positive to
higher dimensional $B$ with polyhedral decompositions, see \cite{PartI},
Definition 1.54. Positivity of $B$
is a necessary condition for $X_0(B,\P,s)$ to appear as the central
fibre of a toric degeneration. All the examples of \S 4 are positive;
this in fact follows from the convexity of reflexive polytopes, and
the positivity condition can be viewed as a type of convexity statement.

Passing back to the general case now, with no restriction on the dimension of
$B$ or the shape of the cells of $\P$, it follows from \cite{PartI}, 
Theorem 3.24, that $\shLS_{X_0(B,\P,s)}$
is a subsheaf of \emph{sets} of a coherent
sheaf we call $\shLS^+_{\pre,X_0(B,\P,s)}$.
This sheaf is a direct sum $\bigoplus_{\omega\in\P\atop\dim\omega=1}
\shN_{\omega}$, where $\shN_{\omega}$ is a line bundle on
the toric stratum of $X_0(B,\P,s)$ corresponding to $\omega$. Furthermore,
as in the two-dimensional normal crossings case, $\shN_{\omega}$ is a
semi-ample line bundle if $B$ is positive.

A section $t\in\Gamma(X_0(B,\P,s),\shLS^+_{\pre,X_0(B,\P,s)})$
which is a section of $\shLS_{X_0(B,\P,s)}$ outside of the zero set $Z$
of $t$ determines a log smooth structure on $X_0(B,\P,s)\setminus Z$.
In particular, if $Z$ does not contain any toric stratum, we are in good shape.
We then obtain a log morphism $X_0(B,\P,s)^{\dagger}\rightarrow 0^{\dagger}$
which is log smooth away from $Z$. We call such a structure a 
\emph{log Calabi-Yau space}.

Let's review: given data
\begin{itemize}
\item $s$ open gluing data;
\item $t\in\Gamma(X_0(B,\P,s),\shLS^+_{\pre,X_0(B,\P,s)})$,
with $t$ a section of $\shLS_{X_0(B,\P,s)}$ over $X_0(B,\P,s)\setminus
Z$ for some set $Z$ which does not contain any toric stratum of $X_0(B,\P,s)$;
\end{itemize}
we obtain $X_0(B,\P,s)^{\dagger}\rightarrow 0^{\dagger}$.

Conversely, we show in \cite{PartI} that if $\X\rightarrow D$ is
a toric degeneration, then $\X_0^{\dagger}\rightarrow 0^{\dagger}$ is
obtained in this way from the dual intersection complex $(B,\P)$
from some choice of data $s$ and $t$. To complete this picture, it
remains to answer

\begin{question}[The reconstruction problem, Version III]
\label{smoothingconjecture}
Suppose $(B,\P)$ is positive.
\begin{enumerate}
\item
What are the possible choices of $s$ and $t$ which
yield log Calabi-Yau spaces?
\item Given $X_0(B,\P,s)^{\dagger}\rightarrow 0^{\dagger}$,
when does it arise as the central fibre of a toric degeneration $\X\rightarrow
D$?
\end{enumerate}
\end{question}

As we have sketched it, this question is
now the refined version of our basic reconstruction problem
Question \ref{reconstruct2}.

\bigskip

The choice of the data $s$ and $t$ determine the moduli of log Calabi-Yau
spaces arising from a given dual intersection
complex. So far we haven't even made the claim that this moduli space is
non-empty, and for general choice of $(B,\P)$, I do not know if this
is the case or not, though it is non-empty if $\dim B=2$ or $3$. However,
one would like a more explicit description of this moduli space in any
event. In general the moduli space is a scheme, but it may be singular
(an example is given in \cite{PartI}, Example 4.28). Some additional hypotheses
are necessary to solve this problem. To motivate the necessary hypothesis,
let's go back to \S 1, where we introduced the notion of simplicity.
We saw that the basic topology of mirror symmetry works only when
the fibration is simple. So maybe we should expect the current construction
to work better when we have simplicity.

There is one technical problem with this: the definition of simplicity
assumes the existence of a torus fibration $f:X\rightarrow B$.
Instead, we want to define simplicity entirely in terms of $B$ itself.
Unfortunately, the solution to this is rather technical, and produces
a definition which is very difficult to absorb (Definition 1.60 of \cite{PartI}). 
Let us just say here that if $B$ is simple in this new sense and
$X(B_0)\rightarrow B_0$ was compactified in a sensible manner
to a topological torus fibration $f:X(B)\rightarrow B$, then $f$ would
be simple in the sense of \S 1, provided that $\dim B\le 3$. In higher
dimensions, this new simplicity does not necessarily imply the old
simplicity; see the forthcoming Ph.D.\ thesis of Helge Ruddat. 
This arises in situations where
orbifold singularities arise in $X(B)$; as is well-known,
such singularities cannot be avoided in higher dimension.

Once we accept this definition, life simplifies a great deal. Extraordinarily,
the a priori very complicated moduli space of log Calabi-Yau spaces with
a given dual intersection complex has a very simple description when
$B$ is simple! One very difficult main result of \cite{PartI},
(Theorem 5.4) is

\begin{theorem}
Given $(B,\P)$ positive and simple, the set of
log Calabi-Yau spaces with dual intersection complex $(B,\P)$,
modulo isomorphism preserving $B$, is $H^1(B,i_*\Lambda\otimes k^{\times})$.
An isomorphism is said to preserve $B$ if it induces the identity on
the dual intersection complex.
\end{theorem}

So the moduli space is an algebraic torus (or a disjoint union of algebraic
tori) of dimension equal to $\dim_k H^1(B,i_*\Lambda\otimes k)$.

Note that this is the expected dimension predicted by the SYZ conjecture.
Indeed, if a smoothing of $X_0(B,\P,s)^{\dagger}$ exists and it
was a topological compactification $X(B)$ of $X_0(B)$,
with a simple torus fibration $f:X(B)\rightarrow B$ extending 
$f_0:X(B_0)\rightarrow B_0$, then $R^{n-1}f_{0*}\RR\cong\Lambda_{\RR}$,
so by simplicity, $R^{n-1}f_*\RR\cong i_*\Lambda_{\RR}$. The discussion
of \S 1 suggests that $\dim H^1(B,R^{n-1}f_*\RR)$ is $h^{1,n-1}$ of the
smoothing, which is of course the dimension of the complex moduli space
of the smoothing.

This argument can be made rigorous by introducing \emph{log differentials}.

\begin{definition} Let $\pi:X^{\dagger}\rightarrow S^{\dagger}$ be
a morphism of logarithmic spaces. A \emph{log derivation}
on $X^{\dagger}$ over $S^{\dagger}$ with values in
an $\O_X$-module $\shE$ is a pair $(\Di,\Dlog)$, where $\Di: \O_X\to
\shE$ is an ordinary derivation of $X/S$ and $\Dlog: \M^\gp_X\to
\shE$ is a homomorphism of abelian sheaves with
$\Dlog\circ\pi^\#=0$; these fulfill the following compatibility
condition
\[
\Di\big(\alpha_X(m)\big)=\alpha_X(m)\cdot \Dlog(m),
\]
for all $m\in\M_X$.

We denote by $\Theta_{X^{\dagger}/S^{\dagger}}$ the sheaf of log derivations
of $X^{\dagger}$ over $S^{\dagger}$ with values in $\O_X$. We set
$\Omega^1_{X^{\dagger}/S^{\dagger}}=\Hom_{\O_X}(\Theta_{X^{\dagger}/
S^{\dagger}},\O_X)$.
\end{definition}

This generalizes the more familiar notion of differentials with
logarithmic poles along a normal crossings divisor.
If $Y\subseteq X$ is a normal crossings divisor inducing
a log structure on $X$, then $\Omega^1_{X^{\dagger}/k}$ is the
sheaf of differentials with logarithmic poles along $Y$, and
$\Omega^1_{Y^{\dagger}/k}$ is the restriction of this sheaf to $Y$.
In general, $\Omega^1_{X^{\dagger}/S^{\dagger}}$ is locally free if
$\pi$ is log smooth. As a result, one can do deformation theory in
the log category for log smooth morphisms (see \cite{F.Kato}). 
This is one of
the principal reasons for introducing log geometry into our picture. 

If $X_0(B,\P,s)^{\dagger}\rightarrow 0^{\dagger}$ is a log Calabi-Yau
space, then the morphism is log smooth off of $Z$. Define
\begin{eqnarray*}
\Theta^p_{X_0(B,\P,s)}&:=&j_*{\bigwedge}^p\Theta_{(X_0(B,\P,s)^{\dagger}
\setminus Z)/0^{\dagger}}\\
\Omega^p_{X_0(B,\P,s)}&:=&j_*{\bigwedge}^p\Omega^1_{(X_0(B,\P,s)^{\dagger}
\setminus Z)/0^{\dagger}}
\end{eqnarray*}
where $j:X_0(B,\P,s)\setminus Z\rightarrow  X_0(B,\P,s)$ is the inclusion.

Then one has

\begin{theorem}
\label{hodgedecomp}
Suppose $(B,\P)$ is positive and simple, and suppose we are given
a log Calabi-Yau space $X_0(B,\P,s)^{\dagger}\rightarrow 0^{\dagger}$ 
which occurs as the central fibre of a toric degeneration $\X\rightarrow
D$ whose general fibre $\X_t$ is non-singular. Then 
for $q=0,1,n-1$ and $n$ with $n=\dim B$,
we have isomorphisms
\begin{eqnarray*}
H^p(B,i_*{\bigwedge}^q\Lambda\otimes k)&\cong&H^p(X_0(B,\P,s),\Theta^q_{X_0(B,\P,s)})
\cong H^p(\X_t,\Theta^q_{\X_t})\\
H^p(B,i_*{\bigwedge}^q\check\Lambda\otimes k)&\cong&H^p(X_0(B,\P,s),\Omega^q_{X_0(B,\P,s)})
\cong H^p(\X_t,\Omega^q_{\X_t})
\end{eqnarray*}
where $\Theta^q_{\X_t}$ and $\Omega^q_{\X_t}$ are the ordinary sheaves
of holomorphic poly-vector fields and holomorphic differentials on a smooth
fibre $\X_t$.
\end{theorem}

The proof of this result, along with a number of other results, appears
in \cite{PartII}.
The result holds for all $q$ when additional hypotheses
are assumed, essentially saying the mirror to $\X_t$ is non-singular.
Note in particular,
since $\Lambda$ and $\check\Lambda$ are interchanged under
discrete Legendre transform, we get the interchange of ordinary
Hodge numbers from this result. In the more general situation where
the Calabi-Yaus arising are singular, one might speculate about
the relationship between these groups, the actual Hodge
numbers and stringy Hodge numbers. These issues are addressed in the
forthcoming Ph.D.\ thesis of Helge Ruddat.

\section{The cone picture and the fan picture}

This section is purely philosophical. In most of our discussion
in \S\S 7 and 8, we focused on the dual intersection complex, and in
particular, focused on the question of constructing a degeneration
from its dual intersection complex. Since
our primary goal was to solve the reconstruction problem 
Question \ref{reconstruct1},
and as the dual intersection complex is related to the complex structure
(Theorem \ref{complextheorem}) it seems natural to focus on the
dual intersection complex. We will see in the next section that this
intuition may not always be correct. So far, the intersection complex
only seemed to arise when talking about mirror symmetry. However,
mirror symmetry instructs us to view both sides of the picture
on the same footing. When we construct a degenerate Calabi-Yau space from a
dual intersection complex, we say we are in the \emph{fan picture},
while if we construct a degenerate Calabi-Yau space from an
intersection complex, we say we are in the \emph{cone picture}.

More precisely, we have seen how given an integral affine manifold
with singularities with toric polyhedral decomposition $(B,\P)$,
then an additional choice of open gluing data $s$ specifies
a space $X_0(B,\P,s)$, along with a sheaf of monoids $\Mbar_{X_0(B,\P,s)}$.
Some additional data may specify a log structure on $X_0(B,\P,s)$ with
this ghost sheaf. The irreducible components of $X_0(B,\P,s)$ are defined
using fans, given by the fan structure of $\P$ at each vertex of 
$\P$. This is why we call this side the fan picture.

On the other hand, given $(B,\P)$ we can also construct 
a \emph{projective} scheme $\check X_0(B,\P,\check s)$ given suitable
gluing data $\check s$. The irreducible components of this scheme
are in one-to-one correspondence with the maximal cells of $\P$;
given such a maximal cell $\sigma$, viewing it as a lattice polytope
in $\RR^n$ determines a projective toric variety, and $\check X_0(B,\P,
\check s)$ is obtained by gluing together these projective toric varieties
using the data $\check s$. This is not quite the same data as occurred
in the fan picture, because we also need to glue the line bundles,
and this is additional data. The reason for calling this side the
cone picture is that each irreducible component can be described
as follows. Given $\sigma\subseteq M_{\RR}$, let $P_{\sigma}=C(\sigma)\cap
(M\oplus\ZZ)$. Then the corresponding projective toric variety is
$\Proj \CC[P_{\sigma}]$, where $\CC[P_{\sigma}]$ is graded using
the projection of $P_{\sigma}$ onto $\ZZ$. Hence the irreducible
components and strata arise from cones over elements of $\P$.

We summarize the duality between the cone and fan pictures:
\medskip
{\scriptsize

\begin{tabular}{|l|l|l|}
\hline
&Fan picture&Cone picture\\
\hline
Gluing data yields & $X_0(B,\P,s)$, $\Mbar_{X_0(B,\P,s)}$&
$\check X_0(B,\P,\check s)$, ample line bundle\\
\hline
$\sigma\in\P$,$\dim\sigma=p$ & An $n-p$-dimensional stratum&
A $p$-dimensional stratum\\
& of
$X_0(B,\P,s)$& of
$\check X_0(B,\P,\check
s)$\\
\hline
$\varphi$ multi-valued&An ample line bundle&A sheaf of monoids\\
convex PL function
&on $X_0(B,\P,s)$
&
$\Mbar_{\check X_0(B,
\P,\check s)}$\\
\hline
$H^p(B,i_*{\bigwedge}^q\Lambda\otimes k)$&
$H^p(X_0(B,\P,s),\Theta^q_{X_0(B,\P,s)})$&
$H^p(\check X_0(B,\P,\check s),\Omega^q_{\check X_0(B,\P,\check s)})$\\
\hline
\end{tabular}
}
\medskip

(Some restrictions may apply to gluing data on both sides in order
for $\varphi$ to yield the desired data.)

In particular, mirror symmetry interchanges discrete information
about the log structure (i.e.\ $\Mbar_{X_0(B,\P,s)}$) and discrete
information about the polarization (i.e the class of the line bundle on
each irreducible component).

\section{Tropical curves}

So far we have seen only the most elementary aspects of mirror
symmetry emerge from this algebro-geometric version of SYZ,
e.g. the interchange of Hodge numbers. However, the real interest
in this approach lies in hints that it will provide a natural
explanation for rational curve counting in mirror symmetry.
If we follow the philosophy of the previous section, we need to
identify structures on affine manifolds with singularities which 
in one of the two pictures has to do with rational curves and in the
other picture has to do with periods. I believe the correct structure
to study is that of tropical curves on affine manifolds with singularities
$B$.
See \cite{Mik},\cite{Sturm} for an introduction to tropical curves in $\RR^n$.
Here, we can take $B$ to be tropical, rather than integral; hence
the name.

\begin{definition}
Let $B$ be a tropical affine manifold with singularities with
discriminant locus $\Delta$. Let
$G$ be a weighted, connected finite graph, with its set of 
vertices and edges denoted by $G^{[0]}$ and $G^{[1]}$
respectively, with weight function $w_{G}:G^{[1]}
\rightarrow\NN\setminus\{0\}$. A parametrized tropical curve in
$B$ is a continuous map $h:G\rightarrow B$ satisfying the following
conditions:
\begin{enumerate}
\item For every edge $E\subseteq G$, $h|_{\Int(E)}$ is
an embedding, $h^{-1}(B_0)$ is dense in $\Int(E)$, and there is 
a section $u\in \Gamma(\Int(E),h^*(i_*\Lambda))$ which is 
tangent to $h(\Int(E))$ at every point of $h(\Int(E))\cap B_0$.
We choose this section to be primitive, i.e.\ not an integral multiple of
another section of $h^*(i_*\Lambda)$.
\item
For every vertex $v\in G^{[0]}$, let $E_1,\ldots,E_m\in G^{[1]}$ be
the edges adjacent to $v$. Let $u_i$ be the section of 
$h^*(i_*\Lambda)|_{\Int(E_i)}$ promised by (1), chosen to point away
from $v$.  This defines
germs $u_i\in h^*(i_*\Lambda)_v=(i_*\Lambda)_{h(v)}$.
\begin{enumerate}
\item
If $h(v)\in B_0$, the following \emph{balancing condition} holds
in
$\Lambda_{h(v)}$:
\[
\sum_{j=1}^m w_{G}(E_j)u_j=0.
\]
\item If $h(v)\not\in B_0$, 
then the following balancing condition is satisfied in
$(i_*\Lambda)_{h(v)}$:
\[
\sum_{j=1}^m w_{G}(E_j)u_j=0
\mod (i_*\check\Lambda)_{h(v)}^{\perp}\cap (i_*\Lambda)_{h(v)}.
\]
The latter group is interpreted as follows. Let $b\in B_0$
be a point near $h(v)$, and identify, via parallel transport
along a path between $h(v)$ and $b$, the groups
$(i_*\Lambda)_{h(v)}$ and $(i_*\check\Lambda)_{h(v)}$ with
local monodromy invariant subgroups of $\Lambda_b$ and
$\check\Lambda_b$ respectively. Then $(i_*\check\Lambda)^{\perp}_{h(v)}$
is a subgroup of $\Lambda_b$, and the intersection makes sense. It is
independent of the choice of $b$ and path.
\end{enumerate}
\end{enumerate}
\end{definition}

So tropical curves behave away from the discriminant locus of $B$
much as the tropical curves of \cite{Mik},\cite{Sturm} do, but they may have legs
terminating on the discriminant locus. As we are interested in
the case that $B$ is compact, we do not want legs which go off
to $\infty$. I warn the reader, however, that this definition is
provisional, and the behaviour in (2) (b) may not be exactly what
we want.

Here we see a tropical elliptic curve, the solid dots being points
of the discriminant locus. The legs terminating at these points
must be in a monodromy invariant direction.

\begin{center}
\includegraphics{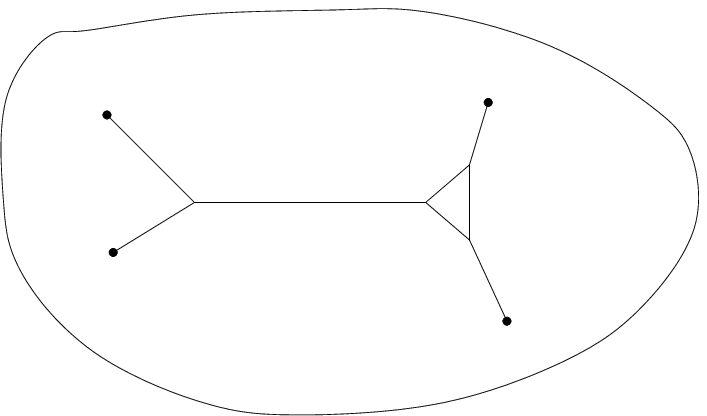}
\end{center}

Now let us connect this to the question of counting curves. In the
situation of a degeneration, $\varphi:\X\rightarrow D$, it is natural
to consider families of maps of curves:
\[
\xymatrix@C=30pt
{\shC\ar[r]^f\ar[d]_{\pi}&\X\ar[d]^{\varphi}\\
D\ar[r]_g&D}
\]
Here $g$ may be a ramified covering, and $\pi$ is a flat morphism
with reduced one-dimensional fibres. In the case of interest,
$f|_{\shC_t}:\shC_t\rightarrow\X_t$ should be a stable map of curves
for $t\not=0$. Let us assume that $\shC_t$ is a non-singular curve for
$t\not=0$. In the logarithmic context, it is then natural to put
the log structure induced by the divisor $\shC_0\subseteq\shC$ on $\shC$,
and so get a diagram
\[
\xymatrix@C=30pt
{\shC^{\dagger}\ar[r]^f\ar[d]_{\pi}&\X^{\dagger}\ar[d]^{\varphi}\\
D^{\dagger}\ar[r]_g&D^{\dagger}}
\]
of log morphisms. Restricting to the central fibre, we obtain a diagram
\[
\xymatrix@C=30pt
{\shC^{\dagger}_0\ar[r]^f\ar[d]_{\pi}&\X^{\dagger}_0\ar[d]^{\varphi}\\
0^{\dagger}\ar[r]_g&0^{\dagger}}
\]
This suggests that we should build up a theory of stable log maps and log
Gromov-Witten invariants. This theory should generalize the theories
developed by Li and Ruan \cite{LR} and Jun Li \cite{Li}.
I will say little about this here, as
this rapidly gets quite technical. There is work in progress of Siebert
on this subject. This point of view has already
been used in \cite{NS} for counting curves in toric varieties, so some
more hints of this approach can be found there. Instead, I wish to sketch
how such a diagram yields a tropical curve.

To do so, consider a
situation where $\pi$ is normal crossings, and the
induced map $\shC^{\dagger}_0\rightarrow\X_0^{\dagger}$ has no
infinitesimal log automorphisms over $0^{\dagger}$. (This is the log
equivalent of the notion of stable map). Let $(B,\P)$ be the dual
intersection complex of the log Calabi-Yau space $\X_0^{\dagger}$.
We can define the \emph{dual intersection graph} of $f:
\shC_0^{\dagger}\rightarrow\X^{\dagger}$, which will be a parameterized
tropical
curve on $B$. I will only do the case here when the image of $f$ is
disjoint from the set $Z\subseteq\X$ of Definition \ref{toricdegen},
(4); otherwise there are some technicalities to worry about.

First we build $G$. Let $C_1,\ldots,C_m$ be the irreducible components
of $\shC_0$. Assume these components are normal for ease of describing
this construction.
Set $G^{[0]}=\{v_1,\ldots,v_m\}$. On the other hand,
$G^{[1]}$ will contain an edge $\overline{v_iv_j}$
joining $v_i$ and $v_j$ whenever
$C_i\cap C_j\not=\emptyset$.

To define $h:G\rightarrow B$, we first describe the image of
each vertex. Let $X_{\sigma_i}$ be the minimal stratum of
$\X_0$ containing $f(C_i)$, where $\sigma_i\in\P$. Let $\eta_i$
be the generic point of $C_i$, $\xi_i=f(\eta_i)$. Then we have
an induced
map $f^{\#}:\M_{\X_0,\xi_i}\rightarrow\M_{\shC_0,\eta_i}$, as $f$ is
a log morphism. This induces a diagram on stalks of 
ghost sheaves
\[
\xymatrix@C=30pt
{\Mbar_{\shC_0,\eta_i}&
\Mbar_{\X_0,\xi_i}
\ar[l]_{\bar f^{\#}}
\\
\Mbar_{0}
\ar[u]^{\bar\pi^{\#}}
&\Mbar_{0}\ar[l]^{\bar g^{\#}}
\ar[u]_{\bar\varphi^{\#}}}
\]
Now $\Mbar_0=\NN$ (see Example \ref{logexamples}, (4))  
and $\Mbar_{\shC_0,\eta_i}=\NN$
since $\pi$ is normal crossings. On the other hand, $\bar\pi^{\#}$
is the identity and if $g$ is a branched cover of degree $d$, then
$g^{\#}$ is multiplication by $d$. By Exercise \ref{ghostexercise},
\[
\Mbar_{\X_0,\xi_i}=\Hom_{monoid}(C(\sigma_i)\cap (M\oplus\ZZ),\NN).
\]
But 
\[
\Hom_{monoid}(\Hom_{monoid}(C(\sigma_i)\cap(M\oplus\ZZ),\NN),\NN)=C(\sigma_i)
\cap(M\oplus\ZZ),
\]
so $\bar f^{\#}$ is determined by an element $(m,r)$ of $C(\sigma_i)
\cap (M\oplus\ZZ)$. Now 
\[
\bar f^{\#}(\bar\varphi^{\#}(1))=\bar f^{\#}(0,1)=\langle (m,r),(0,1)\rangle =r
\]
while
\[
\bar\pi^{\#}(\bar g^{\#}(1))=\bar\pi^{\#}(d)=d.
\]
Thus $r=d$, and $m/d\in\sigma_i$. We define $h(v_i)=m/d$. This is a point
of $\sigma_i$ which is contained in $B$.

If $C_i\cap C_j\not=\emptyset$, there is a minimal stratum $X_{\sigma_{i,j}}$
containing $C_i\cap C_j$. Of course $X_{\sigma_{i,j}}\subseteq X_{\sigma_i}
\cap X_{\sigma_j}$. In particular, $\sigma_{i,j}$ contains $\sigma_i$
and $\sigma_j$.
We take $h(\overline{v_i
v_j})$ to be the straight line joining $h(v_i)$ and $h(v_j)$ inside
$\sigma_{i,j}$. Furthermore, if $\sigma_{i,j}\subseteq M_{\RR}$ is embedded
as a lattice polytope, let $m_{ij}$ be a primitive lattice element
parallel to $m_i-m_j$, and we take $w_{G}(\overline{v_iv_j})$
to be defined by the equation
\[
w_{G}(\overline{v_iv_j})m_{ij}=\#(C_i\cap C_j)(m_i-m_j).
\]

\begin{proposition}
$h$ is a parametrized tropical curve.
\end{proposition}

We do not give a proof here. The case where $X_{\sigma_i}$ is always
an irreducible component of $\X_0$ is essentially covered in \cite{NS}.
Instead, we'll do another extremal case, which exhibits some interesting
features of log geometry.

\begin{example} Suppose a component $C_1$ of $\shC_0$ and all
components $C_2,\ldots,C_t$ intersecting $C_1$ are mapped by $f$
to a zero dimensional stratum $X_{\sigma}$ of $\X_0$. Without loss
of generality we can assume 
\[
\X_0=V(\sigma)=\Spec\CC[\dual{C(\sigma)}\cap(N\oplus\ZZ)]/(z^{(0,1)})
\]
as defined in \S 8. Thus $h$ maps $v_1,\ldots,v_t$ into points
$m_1/d,\ldots,m_t/d\in\sigma$. Let us understand why the balancing
condition holds at $m_1/d$. Let $U\subseteq\shC_0$ be an open neighbourhood
of $C_1$ which only intersects $C_1,\ldots,C_t$, so $h$ is constant
on $U$ as an ordinary morphism (but not as a log morphism). 
Restrict the log structure on $\shC_0$ to $U$. We have
an exact sequence
\[
1\mapright{}\O_U^{\times}\mapright{}\M_U^{\gp}\mapright{p}\Mbar^{\gp}_U
\mapright{} 0.
\]
Taking global sections, we get
\[
1\mapright{}\Gamma(U,\O_U^{\times})\mapright{}\Gamma(U,\M_U^{\gp})
\mapright{p}\Gamma(U,\Mbar_U^{\gp})=
\ZZ^t\mapright{q}\Pic U.
\]
A section $s\in\Gamma(U,\Mbar_U^{\gp})$ defines an $\O_U^{\times}$-torsor
$p^{-1}(s)$, whose class in the Picard group of $U$ is $q(s)$.
It is an easy exercise in log geometry to show that if $s$ is the $i$th
standard basis vector for $\ZZ^t$, then $q(s)=\O_{\shC}(-C_i)|_U$.
Note
$\deg \O_{\shC}(-C_i)|_{C_1}=-\# C_1\cap C_i$ for $i=2,\ldots,t$
and 
$\deg \O_{\shC}(-C_1)|_{C_1}=\sum_{i=2}^t\# C_1\cap C_i$
as $C_1.\shC_0=0$ in $\shC$.

Now observe that $f^{\#}$ acting on the sheaves of monoids induces a diagram
\[
\xymatrix@C=30pt
{\O_{V(\sigma),x}^{\times}\oplus\Mbar^{\gp}_{V(\sigma),x}\ar[r]^{\cong}&
\M^{\gp}_{V(\sigma),x}\ar[r]^{f^{\#}}\ar[d]&\Gamma(U,\M_U^{\gp})\ar[d]^p&\\
N\oplus\ZZ\ar[r]_{\cong}&\Mbar^{\gp}_{V(\sigma),x}\ar[r]_{\bar f^{\#}}&\Gamma(U,
\Mbar_U^{\gp})\ar[r]_{\cong}&\ZZ^t.}
\]
The map $\bar f^{\#}$, by construction, is given by
\[
(n,r)\in N\oplus\ZZ\mapsto (\langle (n,r),(m_i,d)\rangle)_{i=1,\ldots,t}.
\]
On the other hand, in order for $\bar f^{\#}$ to lift to $f^{\#}$, the 
$\O_U^{\times}$ torseur $p^{-1}(\bar f^{\#}(n,r))$ must have a section
for every $(n,r)\in N\oplus\ZZ$, i.e.\ must be trivial in the Picard group.
This implies
\[
\deg\bigotimes_{i=1}^t (\O_{\shC}(-C_i)|_{C_1})^{\otimes \langle (n,r),
(m_i,d)\rangle}=0,
\]
or
\[
\sum_{i=2}^t(\# C_1\cap C_i)(\langle (n,r),(m_i,d)-(m_1,d)\rangle)=0
\]
for all $(n,r)\in N\oplus\ZZ$.
But this is equivalent to
\[
\sum_{i=2}^t (\# C_1\cap C_i) (m_i-m_1)=0,
\]
which is the balancing condition.
\qed
\end{example}

\bigskip

Following the logic of mirror symmetry, this suggests that tropical curves
on the cone side should have to do with periods. It is only recently
that an understanding of this has begun to emerge, and unfortunately,
I do not have space or time to elaborate on this. Let me say that
in \cite{BigPaper}, Siebert and I 
have given a solution 
to Question \ref{smoothingconjecture}, (2), given some hypotheses
on $X_0(B,\P,s)^{\dagger}$, which are implied by simplicity of $B$. 
In this solution, we construct
explicit deformations of a log Calabi-Yau space, order by order. 
Formally, our construction looks somewhat similar to that of Kontsevich
and Soibelman \cite{KS2} for constructing non-Archimedean K3 surfaces
from affine manifolds, and we apply
a key lemma of \cite{KS2}. However, Kontsevich and Soibelman work 
on what we would call the fan side, while we work on the cone side.
This may be surprising given that all of
our discussions involving the strategy of building log Calabi-Yau
spaces was done on the fan side. However, if we take the mirror philosophy
seriously, and we want to see tropical curves appear in a description
of a smoothing, we need to work on the cone side. It turns out to
be extremely natural. In fact,
all tropical rational curves play a role in our
construction. Ultimately, all periods can be calculated in terms of
the data involved in our construction, and in particular, there is a
clear relationship between the period calculation and the existence
of tropical rational curves on $B$.
Once this is fully understood, this will finally give
a firm understanding of a geometric explanation of mirror symmetry.

\end{document}